\documentclass[11pt]{article}

\usepackage{fullpage}
\usepackage{latexsym}
\usepackage{amsfonts}
\usepackage{amssymb}
\usepackage{amsmath}
\usepackage{subcaption}
\usepackage{color}
\usepackage{graphicx, graphics}
\usepackage[colorlinks=true]{hyperref}
\usepackage{epsfig}
\usepackage{xargs}
\usepackage{arydshln}
\newcommand{\sgn}{\mbox{sgn}}
\newcommand{\ud}{\,\mathrm{d}}
\renewcommand\Re{\operatorname{Re}}
\renewcommand\Im{\operatorname{Im}}

\def\Xint#1{\mathchoice
{\XXint\displaystyle\textstyle{#1}}%
{\XXint\textstyle\scriptstyle{#1}}%
{\XXint\scriptstyle\scriptscriptstyle{#1}}%
{\XXint\scriptscriptstyle\scriptscriptstyle{#1}}%
\!\int}
\def\XXint#1#2#3{{\setbox0=\hbox{$#1{#2#3}{\int}$}
\vcenter{\hbox{$#2#3$}}\kern-.5\wd0}}

\def\dashint{\Xint-}

\newcommand{\CC}{\mathbb{C}}

\newcommand{\RR}{\mathbb{R}}

\begin{document}
\title{The time-dependent Schr\"odinger equation with piecewise constant potentials}
\date{\today}


\author{Natalie E. Sheils$^1$ and Bernard Deconinck$^2$ \\
\footnotesize 1. School of Mathematics, University of Minnesota, Minneapolis, MN 55455\\
\footnotesize 2. Department of Applied Mathematics, University of Washington, Seattle, WA 98195-3925 \\
\footnotesize email\textup{: \texttt{nesheils@umn.edu}}
}



\maketitle
%
%
%
%
%

\begin{abstract}
The linear Schr\"odinger equation with piecewise constant potential in one spatial dimension is a well-studied textbook problem.  It is one of only a few solvable models in quantum mechanics and shares many qualitative features with physically important models.  In examples such as ``particle in a box" and tunneling, attention is restricted to the time-independent Schr\"odinger equation.  This paper combines the Unified Transform Method and recent insights for interface problems to present fully explicit solutions for the time-dependent problem.
\end{abstract}


\section{Introduction}

The $N$-particle time-dependent (linear) Schr\"odinger equation is given by

\begin{equation}\label{NS}
i\hbar \frac{\partial \psi}{\partial t}=\left(-\sum_{n=1}^N \frac{p_n^2}{2 m_n}+V(x_1, \ldots, x_N, t) \right)\psi.
\end{equation}
Here $\hbar$ is the reduced Planck constant, $x_j$ denotes the $3$-dimensional coordinate vector of the $j^\textrm{th}$ particle with mass $m_j$, $p_j$ denotes the momentum operator $i \hbar \nabla_{x_j}$ for the $j^{\textrm{th}}$ particle, and $V(x_1, \ldots, x_N,t)$ is the $N$-particle potential.  One can argue that~\eqref{NS} is the most important partial differential equation (PDE) in all of mathematical physics. Standard textbooks such as~\cite{FeynmanLecture2, LandauLifshitz, Merzbacher} rightfully emphasize the solution of~\eqref{NS} in simplified settings, so as to build up the intuition using exact solutions and their properties. Favorite textbook scenarios consider the one-particle case $N=1$ in one (1) spatial dimension with time-independent potential $V(x)$. The linear Schr\"odinger (LS) equation reduces to

\begin{equation}\label{tdls}
i \hbar \frac{\partial \psi}{\partial t}=-\frac{\hbar^2}{2m}\psi_{xx}+V(x)\psi,
\end{equation}
where $m$ is the particle mass. Since $V(x)$ is time independent, separation of variables $\psi(x,t)=\phi(x)T(t)$ leads to

\begin{equation}\label{tils}
T(t)=T_0 e^{-i E t/\hbar}, ~~-\frac{\hbar^2}{2m}\phi''+V(x)\phi=E\phi,
\end{equation}
where the energy $E$ is a (real) separation constant. The second equation above is the one-dimensional one-particle time-independent Schr\"odinger equation. Even at this point, the problem is solvable in closed form in only a few cases, such as the free particle ($V=0$) and the harmonic oscillator ($V=kx^2/2$, $k$ constant)~\cite{FeynmanLecture2, LandauLifshitz, Merzbacher}.

The study of Schr\"odinger equations with piecewise constant potentials is important for a number of reasons. First, to some extent (see below), analytical solutions are available, allowing the development of more physical intuition using scenarios such as the particle in a box, and the piecewise constant potential barrier \cite{Merzbacher}. Piecewise constant potentials also provide the simplest example of a periodic potential, using the Kronig-Penney model~\cite{Merzbacher}. Second, multiple-scale perturbation theory \cite{BenderOrszag, Janowicz, KevorkianCole_MS} shows that the approximation of a complicated $x$-dependent potential using a few constant levels results in accurate leading-order behavior, provided the levels are adequately chosen. This is also evident from the Rayleigh-Ritz characterization of the eigenvalues of \eqref{tils} \cite{LandauLifshitz, Merzbacher}, which depends only on weighted averages of the potential.  As such, the understanding of~\eqref{tdls} with piecewise constant potential is of central importance to the study of quantum mechanics. From a physical point of view, the qualitative features of a potential can often be approximated well using a potential which is pieced together from a number of constant parts~\cite{Griffiths, Merzbacher}.  For instance, although the forces acting between a proton and a neutron are not accurately known on theoretical grounds, it is known that they are short-range forces, \emph{i.e.}, they extend a short distance, then drop to zero quickly.  These forces are well modeled using a piecewise constant potential~\cite{Merzbacher}.

Nonetheless, the solutions that are found in the piecewise constant setting are often restricted to single-mode solutions of~\eqref{tils}, explaining the phenomena of tunneling and trapping~\cite{FeynmanLecture2, LandauLifshitz, Merzbacher}. Solutions of the initial-value problem (IVP) for~\eqref{tdls} are not readily available. The presence of both discrete and continuous spectrum exacerbates the use of straightforward linear superposition. Extensive discussions of this are found in \cite{LevitanSargsjan_SpectralTheory, LevitanSargsjan}, but even there the required superposition result is not immediately found. The goal of this paper is to solve the IVP for~\eqref{tdls} using the Unified Transform Method (UTM) due to Fokas and collaborators~\cite{DeconinckTrogdonVasan, FokasBook, FokasPelloni4}, combined with more recent ideas generalizing the UTM to allow for the explicit solution of interface problems~\cite{Asvestas, DeconinckPelloniSheils, DeconinckSheilsSmith, SheilsDeconinck_PeriodicHeat, SheilsDeconinck_I2I, SheilsDeconinck_LS, SheilsSmith}. In what follows we present {\em explicit, closed-form} solutions of the IVP for~\eqref{tdls} with initial data (and its spatial derivative) that is $L^1(\RR)$ and absolutely continuous on the real line. The solution formulae produced are eminently suitable for asymptotic evaluation, and, if so desired, the location of the discrete and continuous spectrum for the problem may easily be deduced from the solutions presented. It should be noted that the knowledge of these spectra is not required for the construction and evaluation of the solution formulae.

Recently, the UTM has been used to construct explicit closed-form solutions of classical interface problems~\cite{Asvestas, DeconinckPelloniSheils, DeconinckSheilsSmith, SheilsDeconinck_PeriodicHeat, SheilsDeconinck_LS, SheilsSmith}. These are  initial-boundary value problems for partial differential equations (PDEs) for which the solution of an equation in one domain prescribes boundary conditions for the equation in adjacent domains. The standard approach using classical methods to approach such problems is to solve the PDE in each domain, pretending that the boundary values at the domain edges are given. Once solutions in each domain are constructed, the conditions at the interface ({\em e.g.} continuity of the solution and/or some of its derivatives, {\em etc}.) are imposed, resulting in nonlocal equations to be solved for the unknown boundary values at the interfaces. For generic initial conditions, these nonlocal equations are often only solvable numerically. By incorporating the conditions at the interface at an earlier stage, many interface problems can be solved in closed forms~\cite{Asvestas, DeconinckPelloniSheils, DeconinckSheilsSmith, SheilsDeconinck_PeriodicHeat, SheilsDeconinck_LS, SheilsSmith}.

The first problem tackled was that of heat-flow in composite walls or rods \cite{DeconinckPelloniSheils}, or equivalently, diffusion in piecewise homogeneous media~\cite{Asvestas, Mantzavinos}. This was followed by the investigation of an interface problem for the linear free Schr\"odinger equation \cite{SheilsDeconinck_LS}, where we worked with wave functions that are continuous across the interface, but their derivative may experience a jump. We apply the same techniques to the IVP consisting of (\ref{tdls}) with $\psi(x,0)=\psi_0(x)\in L^1(\mathbb{R})$ as well as $\partial_x\psi_0(x)\in L^1(\mathbb{R})$ and both are absolutely continuous on $\mathbb{R}$. We regard this problem as an interface problem with interfaces located at the discontinuities of the potential $V(x)$. The wave function $\psi(x,t)$ and its derivative $\psi_x(x,t)$ are assumed to be continuous across the interfaces. The first condition is a requirement following from the probabilistic interpretation of the wave function, while the second condition follows from integrating the equation across an interface and allowing the length of the integration interval to limit to zero \cite{Merzbacher}. For simplicity, the independent variables occurring in (\ref{tdls}) are rescaled so that, in effect, we may equate $m=1$, $\hbar=1$. Thus in what follows, we consider

\begin{equation}
i \frac{\partial \psi}{\partial t}=-\psi_{xx}+V(x)\psi,\hspace{.8in}x\in\RR,
\end{equation}
where $V(x)$ is a piecewise constant potential. The case where $V(x)$ is a delta function (point singular potential) is covered in~\cite{Rybalko}. We begin our discussions with the case of a single potential jump in Section~\ref{sec:LSp_jump}, where many of the steps of the general method are illustrated in a simple setting. This is followed by the general case of $n$ jumps in Section~\ref{sec:n}. Finally, we apply the results of Section~\ref{sec:n} to the important case of a potential well or barrier in Section~\ref{sec:well}.

In Section~\ref{sec:i2i} we construct a map from the initial conditions to the values of the function and its first spatial derivative evaluated at the $n$ interfaces.  The existence of this map allows one to change the problem at hand from an interface problem to a boundary-value problem (BVP) which allows for an alternative to the approach of finding a closed-form solution to the interface problem.  This was explored previously by the authors for the heat equation in~\cite{SheilsDeconinck_I2I}.

\section{A step potential}\label{sec:LSp_jump}

We wish to solve the classical IVP

\begin{subequations}\label{ls_p}
\begin{align}
&i\psi_t=-\psi_{xx}+\alpha(x)\psi, &-\infty<x<\infty,\\
&\psi(x,0)=\psi_0(x), &-\infty<x<\infty,
\end{align}
\end{subequations}

\noindent where

\begin{equation}
\alpha(x)=\left\{\begin{array}{lcr} \alpha_1,&&x<0,\\\alpha_2,&&x>0, \end{array}\right.
\end{equation}

\begin{figure}
\begin{center}
\def\svgwidth{.7\textwidth}
   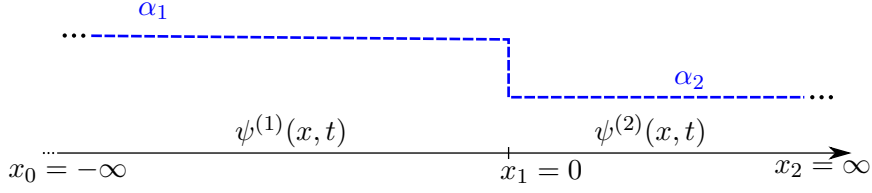 
   \caption{The potential $\alpha(x)$ in the case of one step.   \label{fig:potentialstep}}
   \end{center}
\end{figure}

\noindent $\alpha_1,\alpha_2\in\RR$, and $\lim_{x\to\pm \infty}\psi(x,t)=0$ with $\psi(x,t)$ and its first spatial derivative in $L^1$.  We treat this as an interface problem solved by

\begin{equation}
\psi(x,t)=\left\{ \begin{array}{lcr} \psi^{(1)}(x,t),&&x<0,\\
\psi^{(2)}(x,t),&&x>0, \end{array}\right.
\end{equation}

\noindent where $\psi^{(1)}(x,t)$ and $\psi^{(2)}(x,t)$ solve

\begin{subequations}
\begin{align}
i\psi^{(1)}_t=&-\psi^{(1)}_{xx}+\alpha_1\psi^{(1)}, &x<0,\\
i\psi^{(2)}_t=&-\psi^{(2)}_{xx}+\alpha_2\psi^{(2)}, &x>0,
\end{align}
\end{subequations}

\noindent with initial conditions

\begin{subequations}
\begin{align}
\psi^{(1)}(x,0)=&\psi^{(1)}_0(x), &x<0,\\
 \psi^{(2)}(x,0)=&\psi^{(2)}_0(x), &x>0,
\end{align}
\end{subequations}

\noindent and interface continuity conditions

\begin{subequations}\label{cont2}
\begin{align}
\psi^{(1)}(0,t)=&\psi^{(2)}(0,t), &t>0,\\
\psi^{(1)}_x(0,t)=&\psi^{(2)}_x(0,t), &t>0,
\end{align}
\end{subequations}

\noindent as in Figure~\ref{fig:potentialstep}.

\begin{figure}[htbp]
\begin{center}
\def\svgwidth{4in}
   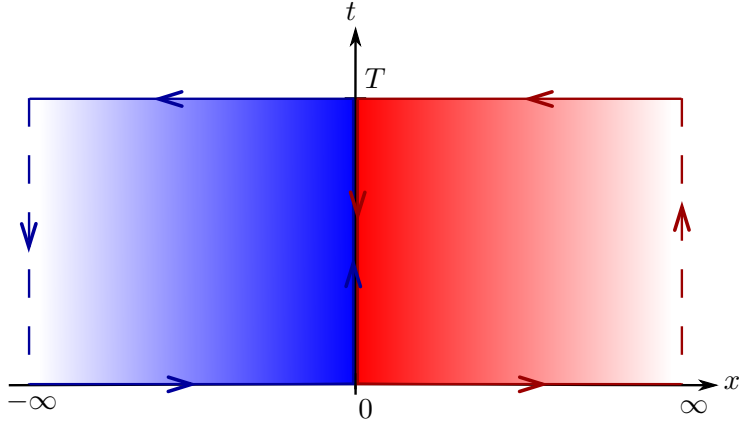 
   \caption{Regions for the application of Green's formula in the case of two semi-infinite domains.   \label{fig:GR_domain2i}}
  \end{center}
\end{figure}

We follow the standard steps in the application of the UTM and begin with the {\em local relations}~\cite{DeconinckTrogdonVasan}:

\begin{subequations}
\begin{align}
(e^{-ikx+\omega_1t}\psi^{(1)})_t=&(e^{-ikx+\omega_1t}(i\psi^{(1)}_x-k\psi^{(1)}))_x, &x<0,\\
(e^{-ikx+\omega_2t}\psi^{(2)})_t=&(e^{-ikx+\omega_2t}(i\psi^{(2)}_x-k\psi^{(2)}))_x, &x>0,
\end{align}
\end{subequations}

\noindent where $\omega_j(k)=i(\alpha_j+k^2)$ for $j=1,2$.  Note that, as is common in the UTM, the $\omega_j$ differ from the standard convention for dispersion relations by a factor of $i$. Thus for dispersive problems $\omega_j$ is purely imaginary. Integrating over the strips $(-\infty,0)\times(0,t)$ and $(0,\infty)\times(0,t)$ respectively (see Figure~\ref{fig:GR_domain2i}), and applying Green's Theorem~\cite{AblowitzFokas}, we have the {\em global relations}

\begin{subequations}
\begin{align}
\int_{-\infty}^0 e^{-ikx+\omega_1t}\psi^{(1)}(x,t)\ud x=&\int_{-\infty}^0e^{-ikx}\psi^{(1)}_0(x)\ud x+\int_0^t e^{\omega_1 s}(i\psi^{(1)}_x(0,s)-k\psi^{(1)}(0,s))\ud s,\\
\int_{0}^\infty e^{-ikx+\omega_2t}\psi^{(2)}(x,t)\ud x=&\int_{0}^\infty e^{-ikx}\psi^{(2)}_0(x)\ud x-\int_0^t e^{\omega_2 s}(i\psi^{(2)}_x(0,s)-k\psi^{(2)}(0,s))\ud s.
\end{align}
\end{subequations}

\noindent Let $\CC^+=\{z\in\CC:\Im(z)\geq0\}$, $\CC^-=\{z\in\CC:\Im(z)\leq0\}$. We define the following:

\begin{align*}
&\hat{\psi}^{(1)}(k,t)=\int_{-\infty}^0 e^{-ikx}\psi^{(1)}(x,t)\ud x, &&x<0,&t>0,&~\Im(k)>0,\\
&\hat{\psi}^{(1)}_0(k)=\int_{-\infty}^0 e^{-ikx}\psi^{(1)}_0(x)\ud x, &&x<0,&&~\Im(k)>0,\\
&\hat{\psi}^{(2)}(k,t)=\int_{0}^\infty e^{-ikx}\psi^{(2)}(x,t)\ud x, &&x>0,&t>0,&~\Im(k)<0,\\
&\hat{\psi}^{(2)}_0(k)=\int_{0}^\infty e^{-ikx}\psi^{(2)}_0(x)\ud x, &&x>0,&&~\Im(k)<0,\\
&g_0(\omega,t)=\int_0^t e^{\omega s} \psi^{(1)}(0,s)\ud s=\int_0^t e^{\omega s} \psi^{(2)}(0,s)\ud s, &&&t>0,&~\omega\in\CC,\\
&g_1(\omega,t)=\int_0^t e^{\omega s} \psi^{(1)}_x(0,s)\ud s=\int_0^t e^{\omega s} \psi^{(2)}_x(0,s)\ud s, &&&t>0,&~\omega\in\CC,
\end{align*}

\noindent where in the last two definitions we have used the continuity conditions \eqref{cont2}. With these definitions the global relations become
\begin{subequations}\label{GR}
\begin{align}
e^{\omega_1t}\hat{\psi}^{(1)}(k,t)=&\hat{\psi}^{(1)}_0(k)+ig_1(\omega_1,t)-kg_0(\omega_1,t), &k\in\CC^+,\label{GR1}\\
e^{\omega_2t}\hat{\psi}^{(2)}(k,t)=&\hat{\psi}^{(2)}_0(k)-ig_1(\omega_2,t)+kg_0(\omega_2,t), &k\in\CC^-. \label{GR2}
\end{align}
\end{subequations}

We wish to transform the global relations so that $g_0(\cdot,t)$ and $g_1(\cdot,t)$ depend on a common argument, $-ik^2$ as was first done in~\cite{Asvestas, Mantzavinos}.  To this end, let
%
%
%
$$
\nu^{(j)}(k)= ik\sqrt{1+\frac{\alpha_j}{k^2}},~~\nu^{(j)}(-k)=-\nu^{(j)}(k),
$$

\noindent which make up a two-sheeted expression with branch points at $\pm i\sqrt{\alpha_j}$ leading to branch cuts in the complex $k$ plane along $[-i\sqrt{\alpha_1},i\sqrt{\alpha_1}]$  and $[-i\sqrt{\alpha_2},i\sqrt{\alpha_2}]$. These cuts are straight-line segments between the endpoints on the real or imaginary axis, depending on the signs of $\alpha_1$ and $\alpha_2$. Note that we have chosen the principal branch, that is, a branch cut along the negative real axis.  Using the transformations $k\rightarrow \nu^{(j)}(\pm k)$, with $j=1$ in~\eqref{GR1} and $j=2$ in~\eqref{GR2},  we have the transformed global relations

\begin{subequations}\label{GRt}
\begin{align}
e^{-ik^2 t}\hat{\psi}^{(1)}\left(\nu^{(1)}(k),t\right)=&\hat{\psi}^{(1)}_0\left(\nu^{(1)}(k)\right)+ig_1(-ik^2,t)-\nu^{(1)}(k)g_0(-ik^2,t),\label{GR1p}\\
e^{-ik^2 t}\hat{\psi}^{(1)}\left(\nu^{(1)}(-k),t\right)=&\hat{\psi}^{(1)}_0\left(\nu^{(1)}(-k)\right)+ig_1(-ik^2,t)-\nu^{(1)}(-k)g_0(-ik^2,t),\label{GR1m}\\
e^{-ik^2t}\hat{\psi}^{(2)}\left(\nu^{(2)}(k),t\right)=&\hat{\psi}^{(2)}_0\left(\nu^{(2)}(k)\right)-ig_1(-ik^2,t)+\nu^{(2)}(k)g_0(-ik^2,t),\label{GR2p}\\
e^{-ik^2t}\hat{\psi}^{(2)}\left(\nu^{(2)}(-k),t\right)=&\hat{\psi}^{(2)}_0\left(\nu^{(2)}(-k)\right)-ig_1(-ik^2,t)+\nu^{(2)}(-k)g_0(-ik^2,t),\label{GR2m}
\end{align}
\end{subequations}
where $\Re(k)\geq0$ in~\eqref{GR1p} and~\eqref{GR2m} and $\Re(k)\leq0$ in~\eqref{GR1m} and~\eqref{GR2p}.

To determine the regions of validity of~\eqref{GRt}, we note that if $\Re(-i \nu^{(j)}(k))$ changes sign then at some point, $-i \nu^{(j)}(k))$ must be purely imaginary.  That is, for some $c\in\RR$:
\begin{align*}
(-i \nu^{(j)}(k)))^2&=-c^2\\
k^2\left(1+\frac{\alpha_j}{k^2}\right)&=-c^2\\
k^2+\alpha_j&=-c^2\\
\Re(k)^2-\Im(k)^2+\alpha_j+2i\Re(k)\Im(k)&=-c^2.
\end{align*}

\noindent Equating real and imaginary parts of the equation, either $\Re(k)=0$ or $\Im(k)=0$.  If $\Re(k)=0$ we have
$$-\Im(k)^2+\alpha_j=-c^2.$$  Thus, $\Re(-i \nu^{(j)}(k))$ can change sign only when $k$ crosses the imaginary axis.  Similarly, if $\Im(k)=0$ we have $\Re(k)^2+\alpha_j=-c^2$ which can only be satisfied if $\alpha_j<0$ and $-\sqrt{-\alpha_j}<\Im(k)<\sqrt{-\alpha_j}$.  In both cases, the sign of $\Re(-i \nu^{(j)}(k))$ is constant for $\Re(k)<0$ and $\Re(k)>0$ (take away $(-\sqrt{|\alpha_j|},\sqrt{|\alpha_j|})$).  By looking at large $k$ asymptotics, we see that $$\sgn(\Re(-i\nu^{(j)}(\pm k)))=\pm\sgn(\Re(k)).$$

\begin{figure}
   \centering
   \def\svgwidth{4in}
   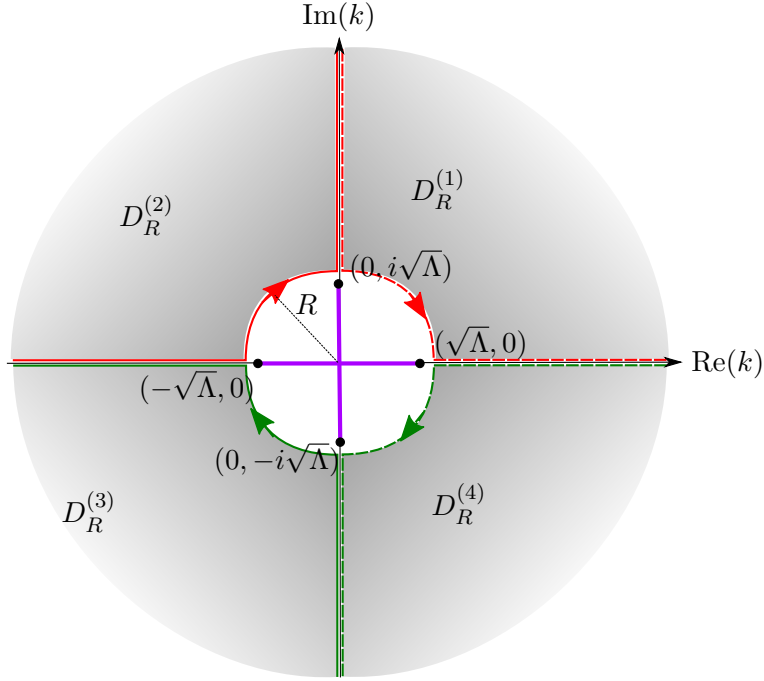
      \caption{The regions $D^{(j)}_R$, $j=1,2,3,4$.
   \label{fig:LSp_DRpm}}
\end{figure}

Inverting the Fourier transform in~\eqref{GR} we have the solution formulae

\begin{subequations}
\begin{align}
\psi^{(1)}(x,t)=&\frac{1}{2\pi}\int_{-\infty}^\infty e^{ikx-\omega_1t}\hat{\psi}^{(1)}_0(k)\ud k+\frac{1}{2\pi}\int_{-\infty}^\infty e^{ikx-\omega_1t} \left(ig_1(\omega_1,t)-kg_0(\omega_1,t)\right)\ud k,\\
\psi^{(2)}(x,t)=&\frac{1}{2\pi}\int_{-\infty}^\infty e^{ikx-\omega_2t}\hat{\psi}^{(2)}_0(k)\ud k-\frac{1}{2\pi}\int_{-\infty}^\infty e^{ikx-\omega_2t}\left(i g_1(\omega_2,t)-kg_0(\omega_2,t)\right)\ud k,
\end{align}
\end{subequations}

\noindent for $x<0$ and $x>0$ respectively.  Examining the second integrals in the formulae above we see it is possible to deform each into the complex plane as follows:
\begin{subequations}\label{reg_dsolns}
\begin{align}
\psi^{(1)}(x,t)=&\frac{1}{2\pi}\int_{-\infty}^\infty e^{ikx-\omega_1t}\hat{\psi}^{(1)}_0(k)\ud k-\frac{1}{2\pi}\int_{\partial D^{(3)}_R} e^{ikx-\omega_1t} \left(ig_1(\omega_1,t)-kg_0(\omega_1,t)\right)\ud k,\\
\psi^{(2)}(x,t)=&\frac{1}{2\pi}\int_{-\infty}^\infty e^{ikx-\omega_2t}\hat{\psi}^{(2)}_0(k)\ud k-\frac{1}{2\pi}\int_{\partial D^{(1)}_R} e^{ikx-\omega_2t}\left(i g_1(\omega_2,t)-kg_0(\omega_2,t)\right)\ud k,
\end{align}
\end{subequations}
where
\begin{equation}\label{DjR}
	D^{(j)}_R = \{ k\in D^{(j)} : |k| > R \},
\end{equation}
with $D^{(j)}$ the $j^{\textrm{th}}$ quadrant of the complex plane.  The regions $D_R^{(j)}$ for $j=1,2,3,4$ are as shown in Figure~\ref{fig:LSp_DRpm} where $\Lambda=\max_l\{ |\alpha_l|\}$ and $R>\sqrt{2\Lambda}$ is a sufficiently large constant.  The reason for integrating around $D_R^{(j)}$ rather than $D^{(j)}$ in~\eqref{reg_dsolns} is to avoid singularities in what follows.

Next we let $k=\nu^{(2)}(\kappa)$ when integrating around $D_R^{(1)}$ and $k=\nu^{(1)}(-\kappa)$ when integrating around $D_R^{(3)}$ so that $g_0(\cdot,t)$ and $g_1(\cdot,t)$ have a common argument and all integrals with unknown terms are integrated around $D_R^{(4)}$.  That is,~\eqref{reg_dsolns} becomes

\begin{subequations}\label{dsolns}
\begin{align}
\begin{split}\label{dsolns1m}
\psi^{(1)}(x,t)=&\frac{1}{2\pi}\int_{-\infty}^\infty e^{ikx-\omega_1t}\hat{\psi}^{(1)}_0(k)\ud k\\
&-\frac{1}{2\pi}\int_{\partial D_{R}^{(4)}} e^{i \nu^{(1)}(-\kappa )x+i\kappa^2t} \left(\frac{i\kappa}{\nu^{(1)}(\kappa)}g_1(-i\kappa^2,t)+\kappa g_0(-i\kappa^2,t)\right)\ud \kappa ,
\end{split}\\
\begin{split}\label{dsolns2p}
\psi^{(2)}(x,t)=&\frac{1}{2\pi}\int_{-\infty}^\infty e^{ikx-\omega_2t}\hat{\psi}^{(2)}_0(k)\ud k\\
&+\frac{1}{2\pi}\int_{\partial D_{R}^{(4)}} e^{i\nu^{(2)}(\kappa)x+i\kappa^2t}\left(\frac{i\kappa}{\nu^{(2)}(\kappa)} g_1(-i\kappa^2,t)-\kappa g_0(-i\kappa^2,t)\right)\ud \kappa.
\end{split}
\end{align}
\end{subequations}
Note that this change of variables maps arcs to arcs but the circular arc of radius $R$ is not mapped exactly to the same circular arc.  However, making another finite contour deformation and using Cauchy's theorem again we may deform to exactly $D_R^{(4)}$.

\begin{figure}[tb]
   \centering
\def\svgwidth{3.5in}
   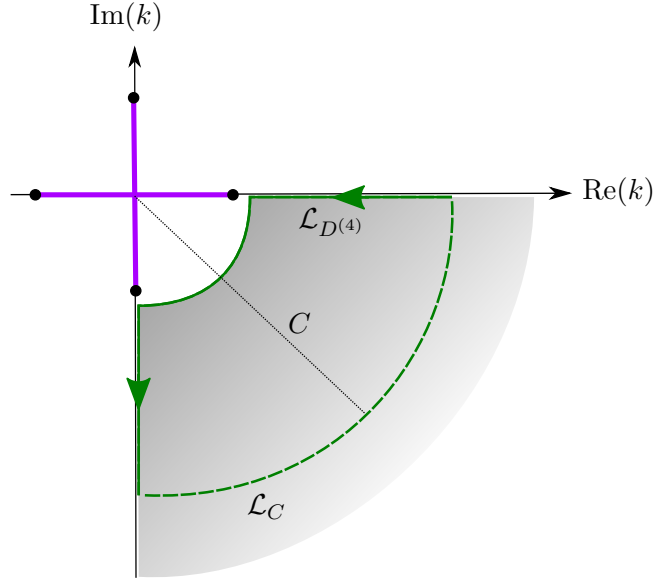
   \caption[The contour $\mathcal{L}^{(4)}$ for the Linear Schr\"odinger equation with potential.]{The contour $\mathcal{L}_{D^{(4)}}$ is shown as a green solid line and the contour $\mathcal{L}_C$ is shown as a green dashed line.  
    \label{LSp_Dpm_close}}
\end{figure}

Using the transformed global relations~\eqref{GR1p} and~\eqref{GR2m} valid in $D^{(4)}$ one solves for $g_0(-i\kappa^2,t)$ and $g_1(-i\kappa^2,t)$.  Noticing that $\nu^{(j)}(-\kappa)=-\nu^{(j)}(\kappa)$ we denote $\nu^{(j)}(\kappa)$ by $\nu^{(j)}$. In the remainder of this section the argument of all $\nu^{(j)}$ is $\kappa$.  Substituting the results for $g_0(-i\kappa^2,t)$ and $g_1(-i\kappa^2,t)$ into~\eqref{dsolns} one finds

\begin{subequations}
\begin{equation}
\label{fullsolns1}
\begin{split}
\psi^{(1)}(x,t)=&\frac{1}{2\pi}\int_{-\infty}^\infty e^{ikx-\omega_1t}\hat{\psi}^{(1)}_0(k)\ud k-\int_{\partial D_{R}^{(4)}} \frac{\kappa(\nu^{(1)}-\nu^{(2)}) }{2\pi\nu^{(1)}(\nu^{(1)}+\nu^{(2)})}e^{-i \nu^{(1)}x+i\kappa^2t} \hat{\psi}^{(1)}_0\left(\nu^{(1)}\right) \ud \kappa \\
&-\int_{\partial D_{R}^{(4)}} \frac{\kappa}{\pi(\nu^{(1)}+\nu^{(2)})} e^{-i \nu^{(1)}x+i\kappa^2t} \hat{\psi}^{(2)}_0\left(-\nu^{(2)}\right) \ud \kappa\\
&+\int_{\partial D_{R}^{(4)}} \frac{\kappa (\nu^{(1)}-\nu^{(2)})}{2\pi\nu^{(1)}(\nu^{(1)}+\nu^{(2)})}\left(\nu^{(2)}-\nu^{(1)}\right)e^{-i \nu^{(1)}x} \hat{\psi}^{(1)}\left(\nu^{(1)},t\right) \ud \kappa \\
&-\int_{\partial D_{R}^{(4)}} \frac{\kappa \nu^{(1)}}{\pi(\nu^{(1)}+\nu^{(2)})} e^{-i \nu^{(1)}x} \hat{\psi}^{(2)}\left(-\nu^{(2)},t\right) \ud \kappa ,
\end{split}
\end{equation}
for $x<0$.  Similarly,

\begin{equation}
\label{fullsolns2}
\begin{split}
\psi^{(2)}(x,t)=&\frac{1}{2\pi}\int_{-\infty}^\infty e^{ikx-\omega_2t}\hat{\psi}^{(2)}_0(k)\ud k-\int_{\partial D_{R}^{(4)}} \frac{\kappa }{\pi(\nu^{(1)}+\nu^{(2)})}\nu^{(2)} e^{i \nu^{(2)}x+i\kappa^2t} \hat{\psi}^{(1)}_0\left(\nu^{(1)}\right) \ud \kappa \\
&+\int_{\partial D_{R}^{(4)}} \frac{\kappa(\nu^{(1)}-\nu^{(2)}) }{2\pi\nu^{(2)}(\nu^{(1)}+\nu^{(2)})} e^{i \nu^{(2)}x+i\kappa^2t} \hat{\psi}^{(2)}_0\left(-\nu^{(2)}\right) \ud \kappa \\
&+\int_{\partial D_{R}^{(4)}} \frac{\kappa }{\pi(\nu^{(1)}+\nu^{(2)})} e^{i\nu^{(2)}x}  \hat{\psi}^{(1)}\left(\nu^{(1)},t\right) \ud \kappa \\
&-\int_{\partial D_{R}^{(4)}} \frac{\kappa(\nu^{(1)}-\nu^{(2)}) }{2\pi\nu^{(2)}(\nu^{(1)}+\nu^{(2)})}e^{i \nu^{(2)}x}  \hat{\psi}^{(2)}\left(-\nu^{(2)},t\right) \ud \kappa ,
\end{split}
\end{equation}
\end{subequations}
for $x>0$.  The first three terms in each of \eqref{fullsolns1} and~\eqref{fullsolns2} depend only on known functions. {The last two terms in~\eqref{fullsolns1} and~\eqref{fullsolns2} are analytic for $\Re(-i\nu^{(1)})=-\Re(i\nu^{(2)})=\Re(\kappa )>0$. 
Note that $\exp(-i\nu^{(1)}x)\hat{\psi}^{(1)}\left(\nu^{(1)},t\right)\to0$ and $\exp(-i\nu^{(1)}x)\hat{\psi}^{(2)}\left(-\nu^{(2)},t\right)\to0$ as $|\kappa|\to\infty$ from within the closure of $D_R^{(4)}$ uniformly in $\kappa$. Thus, by Jordan's Lemma, the integrals of $\exp(-i\nu^{(1)}x)\hat{\psi}^{(1)}\left(\nu^{(1)},t\right)$ and $\exp(-i\nu^{(1)}x)\hat{\psi}^{(2)}\left(-\nu^{(2)},t\right)$ along a closed, bounded curve in the right-half of the complex $\kappa$ plane vanish for $x<0$.  In particular we consider the closed curve  $\mathcal{L}^{(4)}=\mathcal{L}_{D^{(4)}}\cup\mathcal{L}^{(4)}_C$ where $\mathcal{L}_{D^{(4)}}=\partial D_R^{(4)} \cap \{\kappa: |\kappa|<C\}$ and $\mathcal{L}^{(4)}_C=\{\kappa\in D_R^{(4)}: |\kappa|=C\}$, see Figure~\ref{LSp_Dpm_close}.} \label{JordanCauchyArgument}

Since the integral along $\mathcal{L}_C$ vanishes for large $C$, the fourth and fifth integrals on the right-hand side of \eqref{fullsolns1} must vanish since the contour $\mathcal{L}_{D^{(4)}}$  becomes $\partial D^{(4)}$ as $C\to\infty$.  For the final two integrals in Equation~\eqref{fullsolns2} we use the fact that for $x>0$ the integrals of $\exp(i\nu^{(2)}x)\hat{\psi}^{(1)}\left(\nu^{(1)},t\right)$ and $\exp(i\nu^{(2)}x)\hat{\psi}^{(2)}\left(-\nu^{(2)},t\right)$ along a closed, bounded curve in the right-half of the complex $\kappa$ plane vanish. Thus, we have an explicit representation for $\psi^{(1)}(x,t)$ in terms of only initial conditions:
\begin{subequations}
\begin{equation}
\label{totalsoln1}
\begin{split}
\psi^{(1)}(x,t)=&\frac{1}{2\pi}\int_{-\infty}^\infty e^{ikx-\omega_1t}\hat{\psi}^{(1)}_0(k)\ud k\\
&-\int_{\partial D_{R}^{(4)}} \frac{\kappa(\nu^{(1)}-\nu^{(2)}) }{2\pi\nu^{(1)}(\nu^{(1)}+\nu^{(2)})}e^{-i \nu^{(1)}x+i\kappa^2t} \hat{\psi}^{(1)}_0\left(\nu^{(1)}\right) \ud \kappa \\
&-\int_{\partial D_{R}^{(4)}} \frac{\kappa}{\pi(\nu^{(1)}+\nu^{(2)})} e^{-i \nu^{(1)}x+i\kappa^2t} \hat{\psi}^{(2)}_0\left(-\nu^{(2)}\right) \ud \kappa,
\end{split}
\end{equation}

\noindent for $x<0$, and

\begin{equation}
\label{totalsoln2}
\begin{split}
\psi^{(2)}(x,t)=&\frac{1}{2\pi}\int_{-\infty}^\infty e^{ikx-\omega_2t}\hat{\psi}^{(2)}_0(k)\ud k\\
&-\int_{\partial D_{R}^{(4)}} \frac{\kappa \nu^{(2)} }{\pi(\nu^{(1)}+\nu^{(2)})} e^{i \nu^{(2)}x+i\kappa^2t} \hat{\psi}^{(1)}_0\left(\nu^{(1)}\right) \ud \kappa \\
&+\int_{\partial D_{R}^{(4)}} \frac{\kappa(\nu^{(1)}-\nu^{(2)}) }{2\pi\nu^{(2)}(\nu^{(1)}+\nu^{(2)})} e^{i \nu^{(2)}x+i\kappa^2t} \hat{\psi}^{(2)}_0\left(-\nu^{(2)}\right) \ud \kappa ,
\end{split}
\end{equation}
\end{subequations}
\noindent for $x>0$. Note that the denominators in~\eqref{totalsoln1} and~\eqref{totalsoln2} are zero at the branch points $\kappa =\pm i\sqrt{\alpha_j}$.   However, these points are avoided by integrating over the boundary of $D_{R}^{(4)}$.

\begin{figure}
   \centering
   \def\svgwidth{.75\textwidth}
   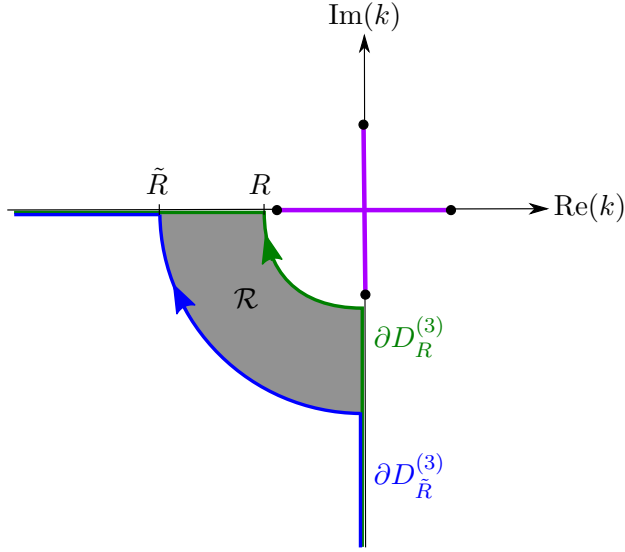
      \caption{The contours $\partial D^{(3)}_R$ and  $\partial D^{(3)}_{\tilde{R}}$ and the region $\mathcal{R}$.
   \label{fig:LSp_DR_noR}}
\end{figure}

The expressions \eqref{totalsoln1} and \eqref{totalsoln2} provide fully explicit solutions for the IVP~\eqref{ls_p}. They are written in a form containing more familiar exponents by  letting $\kappa =ik\sqrt{1+\alpha_1/k^2}$ in the second and third integrals of~\eqref{totalsoln1} and $\kappa=-ik\sqrt{1+{\alpha_2/k^2}}$ in the second and third integrals of~\eqref{totalsoln2}. Then
\begin{subequations}\label{t_totalsolns}
\begin{equation}
\label{t_totalsoln1}
\begin{split}
\psi^{(1)}(x,t)=&\frac{1}{2\pi}\int_{-\infty}^\infty e^{ikx-\omega_1t}\hat{\psi}^{(1)}_0(k)\ud k\\
&-\int_{\partial D_R^{(3)}} \frac{1-\sqrt{1+\frac{\alpha_1-\alpha_2}{k^2}}}{2\pi\left(1+\sqrt{1+\frac{\alpha_1-\alpha_2}{k^2}}\right)} e^{i k x-\omega_1 t} \hat{\psi}^{(1)}_0\left(-k\right) \ud k\\
&-\int_{\partial D_R^{(3)}} \frac{1}{\pi\left(1+\sqrt{1+\frac{\alpha_1-\alpha_2}{k^2}}\right)} e^{i k x-\omega_1t} \hat{\psi}^{(2)}_0\left(k\sqrt{1+\frac{\alpha_1-\alpha_2}{k^2}}\right) \ud k,
\end{split}
\end{equation}

\noindent for $x<0$, and

\begin{equation}
\label{t_totalsoln2}
\begin{split}
\psi^{(2)}(x,t)=&\frac{1}{2\pi}\int_{-\infty}^\infty e^{ikx-\omega_2t}\hat{\psi}^{(2)}_0(k)\ud k\\
&+\int_{\partial D_R^{(1)}} \frac{1}{\pi\left(1+\sqrt{1-\frac{\alpha_1-\alpha_2}{k^2}}\right)} e^{i k x-\omega_2t} \hat{\psi}^{(2)}_0\left(k\sqrt{1-\frac{\alpha_1-\alpha_2}{k^2}}\right) \ud k\\
&+\int_{\partial D_R^{(1)}} \frac{1-\sqrt{1-\frac{\alpha_1-\alpha_2}{k^2}}}{2\pi\left(1+\sqrt{1-\frac{\alpha_1-\alpha_2}{k^2}}\right)} e^{i k x-\omega_2t} \hat{\psi}^{(1)}_0\left(-k\right) \ud k,
\end{split}
\end{equation}
\end{subequations}
for $x>0$.  It appears that our solution depends on an extra parameter $R$. However, observe that $\int_{\partial D^{(3)}_R} \cdot \ud k=\int_{\partial D^{(3)}_{\tilde{R}}} \cdot \ud k+\oint_{\mathcal{R}}$ where $\partial D_R^{(3)}$, $\partial D^{(3)}_{\tilde{R}}$ and $\mathcal{R}$ are as in Figure~\ref{fig:LSp_DR_noR}.  Since the integrands in~\eqref{t_totalsoln1} are analytic in $D_R^{(3)}$ (and therefore $\mathcal{R}$), $\oint_\mathcal{R}\cdot\ud k=0$.  Hence, $\int_{\partial D_R} \cdot \ud k=\int_{\partial D_{\tilde{R}}} \cdot \ud k$ for any $R>\Lambda$, and our solution is independent of the value of $R$ chosen. The same argument is true for~\eqref{t_totalsoln2}.

\begin{figure}
\begin{center}
\begin{subfigure}[b]{.2\textwidth}
\centering
   \def\svgwidth{\textwidth}
   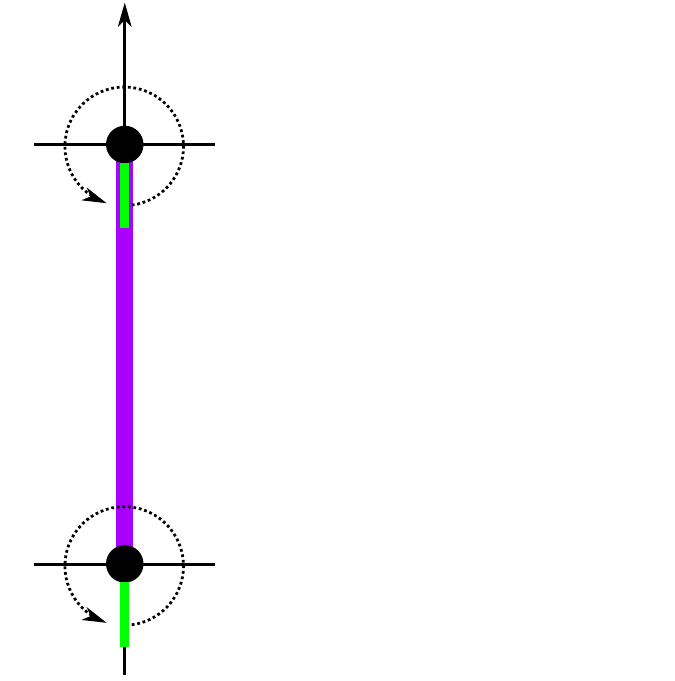
      \caption{\label{fig:LSp_cut_ag0}}
\end{subfigure}
\hfill
\begin{subfigure}[b]{.2\textwidth}
\centering
   \def\svgwidth{\textwidth}
   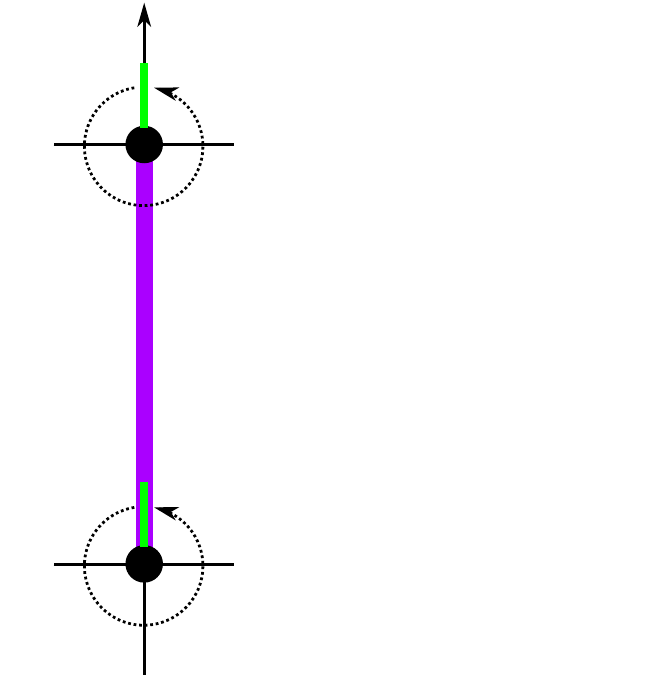
      \caption{\label{fig:LSp_cut_ag0_2}}
\end{subfigure}
\hfill
 \begin{subfigure}[b]{.45\textwidth}
   \centering
   \def\svgwidth{\textwidth}
   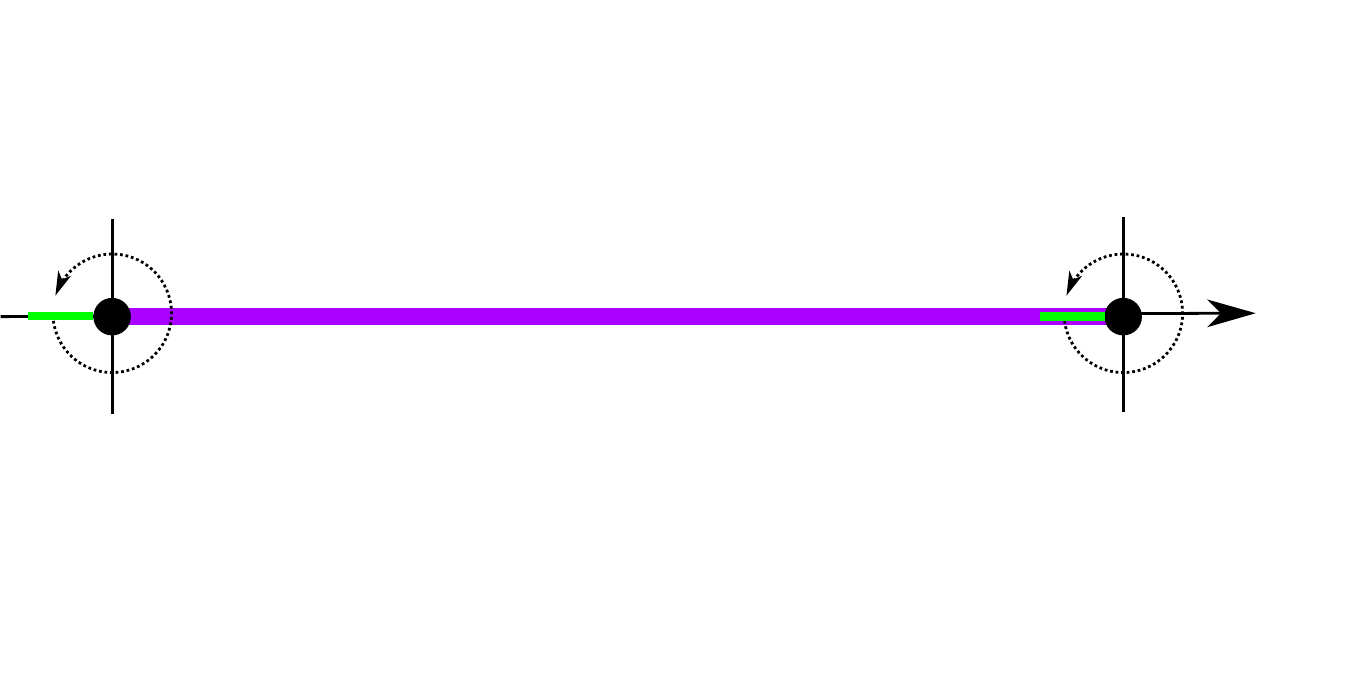
      \caption{ \label{fig:LSp_cut_al0}}
   \end{subfigure}
   \caption{Branch cuts for $\sqrt{1+\frac{a}{k^2}}$ and the local parameterizations around the branch points. In (a), $a>0$ and the local parameterization around the branch points $\pm i\sqrt{a}$ with $-\pi/2<\theta_1,\theta_2\leq 3\pi/2$. In (b), $a>0$ and the local parameterization around the branch points $\pm i\sqrt{a}$ with $-3\pi/2<\theta_3,\theta_4\leq \pi/2$ .  In (c) $a<0$ and the local parameterization around the branch points $\pm \sqrt{-a}$ with $-\pi<\theta_5,\theta_6\leq \pi$.
   \label{fig:LSp_cuts}}
   \end{center}
\end{figure}

Before we begin the asymptotic analysis, it is useful to deform the contours in~\eqref{t_totalsolns} back to the real line.  We examine the branch cut introduced in~\eqref{t_totalsolns} of the form $\sqrt{1+\frac{a}{k^2}}$.  In~\eqref{t_totalsolns} $a=\alpha_2-\alpha_1$ but it may be different in later sections.  If $a>0$ the branch points are at $\pm i\sqrt{a}$.  We fix the branch cut to be on the finite imaginary axis running from $-i\sqrt{a}$ to $i\sqrt{a}$ by defining the local polar coordinates
\begin{align*}
k-i\sqrt{a}=r_1e^{i\theta_1},\\
k+i\sqrt{a}=r_2e^{i\theta_2},
\end{align*}
where $-\pi/2<\theta_1,\theta_2\leq 3\pi/2$ as in Figure~\ref{fig:LSp_cut_ag0} or  $-3\pi/2<\theta_3,\theta_4\leq \pi/2$ as in Figure~\ref{fig:LSp_cut_ag0_2}.  Similarly, if $a<0$, $\pm \sqrt{-a}$ are the branch points.  We fix the branch cut to be on the finite real axis running from $-\sqrt{-a}$ to $\sqrt{-a}$ by defining the local polar coordinates
\begin{align*}
k+\sqrt{-a}=r_5e^{i\theta_5},\\
k-\sqrt{-a}=r_6e^{i\theta_6},
\end{align*}
where $-\pi<\theta_5,\theta_6\leq \pi$ as in Figure~\ref{fig:LSp_cut_al0}.

\begin{figure}
   \centering
   \def\svgwidth{4.5in}
   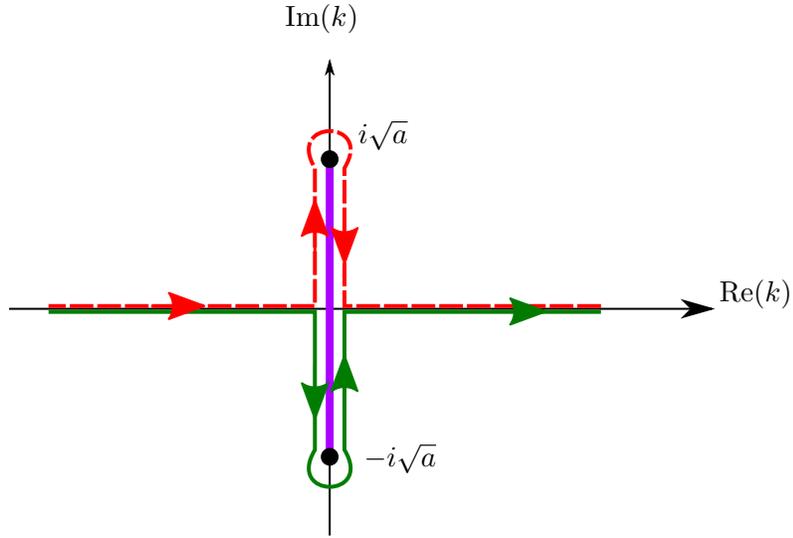
      \caption{The deformations of $\partial D^{(1)}_R$ (as a red dashed line) and $\partial D^{(3)}_R$ (as a green solid line) to the real line when the branch cut is on the imaginary axis.
   \label{fig:LSp_keyhole_ag0}}
\end{figure}

\begin{figure}
   \centering
   \def\svgwidth{4.5in}
   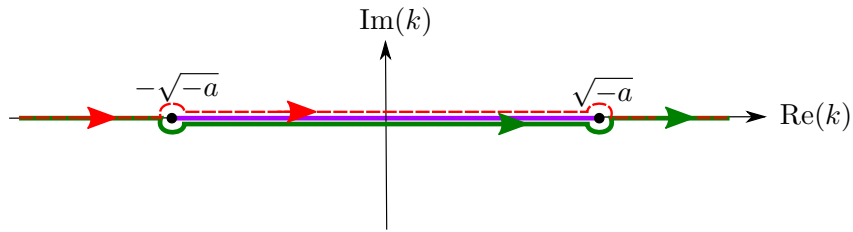
      \caption{The deformations of $\partial D^{(1)}_R$ (as a red dashed line) and $\partial D^{(3)}_R$ (as a solid green line) to the real line for the case when the branch cut is on the real axis.
   \label{fig:LSp_keyhole_al0}}
\end{figure}

If the branch cut is on the imaginary axis then deforming $\partial D_R^{(1)}$ to the real axis and using the local parameterization $-\pi/2<\theta_1,\theta_2\leq 3\pi/2$ as in Figure~\ref{fig:LSp_cut_ag0}, one finds
\begin{equation}\label{1_imag}
\begin{split}
\int_{\partial D_R^{(1)}} f(k) \ud k=&\dashint_{-\infty}^\infty f(k) \ud k+i \int_0^{\sqrt{a}} f(r e^{3\pi i/2}+i\sqrt{a}) \ud r\\
&-\lim_{\epsilon\to 0} i\epsilon \int_{-\pi/2}^{3\pi/2} f(\epsilon e^{i\theta}+i\sqrt{a})e^{i\theta}\ud \theta-i \int_0^{\sqrt{a}} f(r e^{-\pi i/2}+i\sqrt{a}) \ud r,
\end{split}
\end{equation}
as in the dashed red line in Figure~\ref{fig:LSp_keyhole_ag0} where $\dashint_{-\infty}^\infty\cdot \ud k$ is a Cauchy Principal Value integral.  The third integral in~\eqref{1_imag} can be contracted to a zero radius using integration by parts~\cite[p.\ 128]{Miller}.  Deforming $\partial D_R^{(3)}$ to the real axis when the branch cut is on the imaginary axis requires the local parameterization with $-3\pi/2<\theta_3,\theta_4\leq \pi/2$ as in~\ref{fig:LSp_cut_ag0_2}.  Then
\begin{equation}\label{3_imag}
\begin{split}
\int_{\partial D_R^{(3)}} f(k) \ud k=&-\dashint_{-\infty}^\infty f(k) \ud k+i \int_{-\sqrt{a}}^{0} f(r e^{\pi i/2}-i\sqrt{a}) \ud r\\
&-\lim_{\epsilon\to 0} i\epsilon \int_{-3\pi/2}^{\pi/2} f(\epsilon e^{i\theta}-i\sqrt{a})e^{i\theta}\ud \theta-i \int_{-\sqrt{a}}^{0} f(r e^{-3\pi i/2}-i\sqrt{a}) \ud r,
\end{split}
\end{equation}
as in the solid green line in Figure~\ref{fig:LSp_keyhole_ag0}.  Again,  the third integral in~\eqref{3_imag} can be contracted to a zero radius using integration by parts~\cite[p.\ 128]{Miller}. If the branch cut is on the real axis then deforming $\partial D_R^{(1)}$ and $\partial D_R^{(3)}$ to the real axis one finds
\begin{equation}\label{1_real}
\int_{\partial D_R^{(1)}} f(k) \ud k=\dashint_{-\infty}^\infty f(k) \ud k,
\end{equation} and
\begin{equation}\label{3_real}
\int_{\partial D_R^{(3)}} f(k) \ud k=-\dashint_{-\infty}^\infty f(k) \ud k,
\end{equation}
as in Figure~\ref{fig:LSp_keyhole_al0}.

In what follows we consider $\alpha_2>\alpha_1$.  Then, $\partial D_R^{(3)}$ in~\eqref{t_totalsoln1} can be deformed as in~\eqref{3_real} and $\partial D_R^{(1)}$ in~\eqref{t_totalsoln2} can be deformed as in~\eqref{1_imag}.

\begin{equation}
\label{R_totalsoln1}
\begin{split}
\psi^{(1)}(x,t)=&\frac{1}{2\pi}\int_{-\infty}^\infty e^{ikx-\omega_1t}\hat{\psi}^{(1)}_0(k)\ud k+\frac{1}{2\pi}\dashint_{-\infty}^\infty a^{(1)}(k) e^{i k x-\omega_1 t} \ud k,
\end{split}
\end{equation}
for $x<0$, where
$$a^{(1)}(k)=\frac{1}{1+\sqrt{1+\frac{\alpha_1-\alpha_2}{k^2}}} \left(\left(1-\sqrt{1+\frac{\alpha_1-\alpha_2}{k^2}}\right) \hat{\psi}^{(1)}_0\left(-k\right)+2\hat{\psi}^{(2)}_0\left(k\sqrt{1+\frac{\alpha_1-\alpha_2}{k^2}}\right)\right),$$ and

\begin{equation}
\label{R_totalsoln2}
\begin{split}
\psi^{(2)}(x,t)=&\frac{1}{2\pi}\int_{-\infty}^\infty e^{ikx-\omega_2t}\hat{\psi}^{(2)}_0(k)\ud k+\frac{1}{2\pi}\dashint_{-\infty}^\infty a^{(2)}(k) e^{i k x-\omega_2t} \ud k\\
&+i\int_0^{\sqrt{\alpha_2-\alpha_1}} \Big(a^{(2)}(re^{3\pi i/2}+i\sqrt{\alpha_2-\alpha_1})-a^{(2)}(re^{-\pi i/2}+i\sqrt{\alpha_2-\alpha_1})\Big)\\
&~~~~~~~~~~~~~~~~~~~~~~~\times e^{(r-\sqrt{\alpha_2-\alpha_1})x+it(r^2-\alpha_1-2r\sqrt{\alpha_2-\alpha_1})}\ud r,
\end{split}
\end{equation}
for $x>0$, where
$$a^{(2)}(k)= \frac{1}{1+\sqrt{1-\frac{\alpha_1-\alpha_2}{k^2}}}\left( (1-\sqrt{1-\frac{\alpha_1-\alpha_2}{k^2}}) \hat{\psi}^{(1)}_0\left(-k\right)+2\hat{\psi}^{(2)}_0\left(k\sqrt{1-\frac{\alpha_1-\alpha_2}{k^2}}\right)\right).$$

At this point, we are ready to use asymptotic analysis. The large-time leading-order behavior of \eqref{ls_p} with initial conditions which decay sufficiently fast at $\pm \infty$ is easily obtained using the Method of Stationary Phase~\cite{BenderOrszag}.  Notice that the third integrand of~\eqref{R_totalsoln2} is decaying for $x$ large and positive.  Thus, this integral does not contribute using the Method of Stationary Phase.  We choose $x/t=\gamma_1<0$ for $x<0$ and $x/t=\gamma_2>0$ for $x>0$. We obtain
\begin{equation*}
\begin{split}
&\psi^{(1)}\sim \frac{e^{i \left(\frac{\gamma_1^2}{4}-\alpha_1\right)t-\frac{i\pi}{4}}}{2\sqrt{\pi t}}\!\left(\hat{\psi}_0^{(1)}\left(\frac{\gamma_1}{2}\right)\!+\!\frac{\left(1\!-\!\sqrt{1\!+\!\frac{4(\alpha_1-\alpha_2)}{\gamma_1^2}}\right)\hat{\psi}_0^{(1)}\left(\frac{-\gamma_1}{2}\right)\!+\!2\hat{\psi}_0^{(2)}\left(\frac{\gamma_1}{2}\sqrt{1\!+\!\frac{4(\alpha_1-\alpha_2)}{\gamma_1^2}}\right)}{1+ \sqrt{1\!+\!\frac{4(\alpha_1-\alpha_2)}{\gamma_1^2} }} \right),
\end{split}
\end{equation*}
and

\begin{equation*}
\begin{split}
\psi^{(2)}\sim& \frac{e^{i \left(\frac{\gamma_2^2}{4}-\alpha_2\right)t-\frac{i\pi}{4}}} {2\sqrt{\pi t}}\! \left( \hat{\psi}_0^{(2)}\left(\frac{\gamma_2}{2}\right)\!+\!\frac{\left(1\!-\!\sqrt{1\!-\!\frac{4(\alpha_1-\alpha_2)}{\gamma_2^2}}\right)\hat{\psi}_0^{(1)}\left(\frac{-\gamma_2}{2}\right)+2\hat{\psi}_0^{(2)}\left(\frac{\gamma_2}{2}\sqrt{1\!-\!\frac{4(\alpha_1-\alpha_2)}{\gamma_2^2}}\right)}{1+ \sqrt{1\!-\!\frac{4(\alpha_1-\alpha_2)}{\gamma_2^2} }} \right).
\end{split}
\end{equation*}

\begin{figure} 
   \centering
   \def\svgwidth{.75\textwidth}
   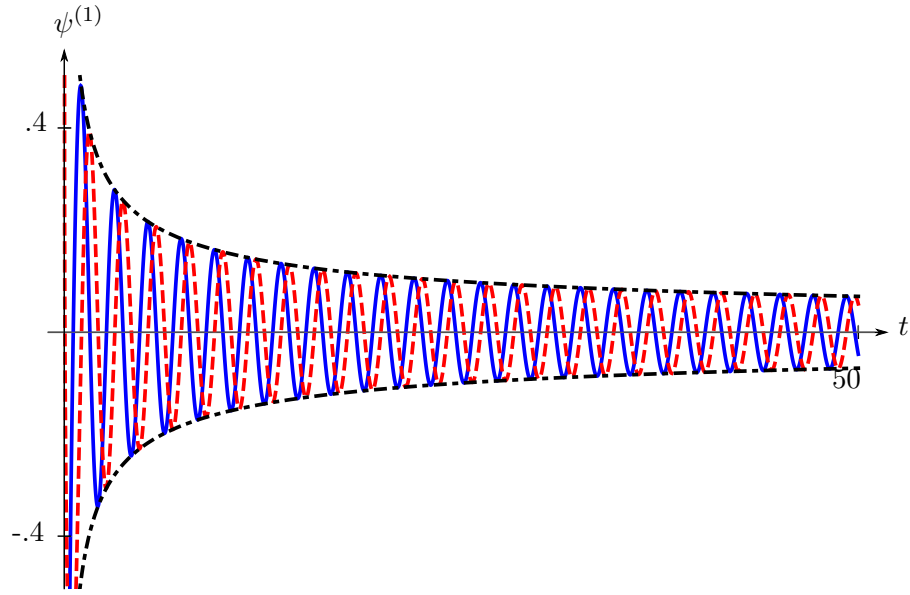
   \caption{The real (red dashed) and imaginary (blue solid) parts of the leading order behavior as $t\to\infty$ of $\psi^{(1)}$ along rays of $x/t=-4$ with $\psi_0(x)=e^{-x^2}$, $\alpha_1=1$, and $\alpha_2=2$.}
   \label{fig:LOB2_1}
\end{figure}

\begin{figure} 
   \centering
   \def\svgwidth{.75\textwidth}
   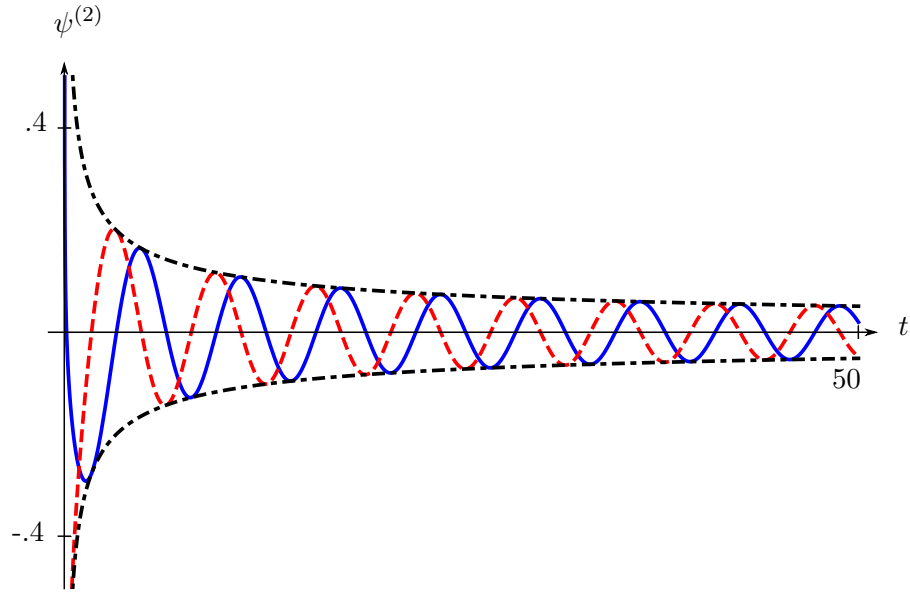
   \caption{The real (red dashed) and imaginary (blue solid) parts of the leading order behavior as $t\to\infty$ of $\psi^{(2)}$ along rays of $x/t=2$ with $\psi_0(x)=e^{-x^2}$, $\alpha_1=1$, and $\alpha_2=2$.}
   \label{fig:LOB2_2}
\end{figure}

The oscillations that are expected as a consequence of dispersion are contained in $\exp(it(\gamma_j^2/4-\alpha_j))$.  In Figures~\ref{fig:LOB2_1} and~\ref{fig:LOB2_2}  the envelopes of the solutions are plotted in black as a dot-dashed line. The real part of the solution (plotted as a solid line in blue) and the imaginary part of the solution (plotted as a dashed line in red) are centered around the $t$-axis.  Using the Method of Stationary Phase one must look in directions of constant $x/t$.  In Figure~\ref{fig:LOB2_1} we consider solutions for $x/t=-4$ and in Figure~\ref{fig:LOB2_2} we have solutions with $x/t=2$. In both figures $\alpha_1=1, \alpha_2=2$ and $\psi_0(x)=e^{-x^2}$.

The Method of Stationary Phase is not useful for considering the nature of solutions near the barrier at $x=0$, since requiring $t$ to be large implies that $x$ is large if $x/t$ is to be constant. In order to evaluate the solution formulae numerically near the interface one could use techniques presented in~\cite{Levin, TrogdonThesis, Trogdon}.  It may also be possible to use asymptotic techniques similar to those in~\cite{BiondiniTrogdon}.


Notice that when $\alpha_1=\alpha_2=0$ the problem reduces to the IVP for the linear Schr\"odinger equation on the whole line. It is easily seen that the solutions~\eqref{totalsoln1} and~\eqref{totalsoln2} reduce to the solution of the problem found using Fourier transforms split into the appropriate domains for the free particle problem.


\section{$n$ potential jumps}\label{sec:n}

We wish to solve the classical problem

\begin{subequations}\label{ls_p_n}
\begin{align}
&i\psi_t=-\psi_{xx}+\alpha(x)\psi, &-\infty<x<\infty,\\
&\psi(x,0)=\psi_0(x), &-\infty<x<\infty,
\end{align}
\end{subequations}

\noindent with

$$
\alpha(x)=\left\{\begin{array}{lrr} \alpha_1,&&x<x_1,\\\alpha_2,&&x_1<x<x_2\\ \vdots&&\\ \alpha_n,&&x_{n-1}<x<x_n,\\ \alpha_{n+1},&&x>x_{n}, \end{array}\right.
$$

\noindent and $\lim_{|x|\to\infty}\psi(x,t)=0$.  Further, recall, $\psi(x,t)$ and its first spatial derivative are both in $L^{1}.$  We repeat the same steps as in the previous section, but now for an arbitrary number $n$ of constant levels of the potential $\alpha(x)$. As a consequence, the formulae obtained are significantly more involved, but no less explicit. The experience gained from the previous section provides the insight necessary to proceed with the general case presented here.

\begin{figure}
\begin{center}
\def\svgwidth{\textwidth}
   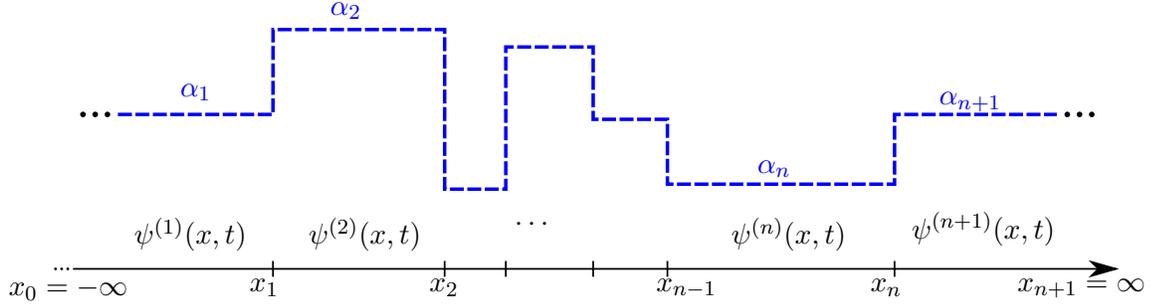 
   \caption{A cartoon of the potential $\alpha(x)$ in the case of $n$ interfaces.   \label{fig:potentialcartoon}}
   \end{center}
\end{figure}

We treat the problem~\eqref{ls_p_n} as an interface problem solved by

\begin{equation}
\psi(x,t)=
\left\{ \begin{array}{lcrrr}
\psi^{(1)}(x,t),&&&&x<x_1,\\
\psi^{(2)}(x,t),&&&x_1<\hspace{-.25cm}&x<x_2,\\
\vdots&&&&\\
\psi^{(n)}(x,t),&&&x_{n-1}<\hspace{-.25cm}&x<x_n,\\
\psi^{(n+1)}(x,t),&&&&x>x_n,
\end{array}\right.
\end{equation}

\noindent which solve the $n+1$ IVPs

\begin{subequations}\label{ivp1}
\begin{align}
i\psi^{(j)}_t&=-\psi^{(j)}_{xx}+\alpha_j \psi^{(j)},\\
\psi^{(j)}(x,0)&=\psi_0^{(j)}(x),
\end{align}
\end{subequations}

\noindent for $x_{j-1}<x<x_j$, with $x_0=-\infty$ and $x_{n+1}=\infty$, $j=1, \ldots, n+1$.  The solutions of the IVPs~\eqref{ivp1} are coupled by the interface conditions

\begin{align*}
\psi^{(j)}(x_j,t)=&\psi^{(j+1)}(x_j,t), &t>0,\\
\psi^{(j)}_x(x_j,t)=&\psi^{(j+1)}_x(x_j,t), &t>0,
\end{align*}

\noindent for $1\leq j\leq n$ as in Figure~\ref{fig:potentialcartoon}.

We begin with the $n+1$ local relations

\begin{equation}\label{n_lr}
(e^{-ikx+\omega_jt}\psi^{(j)})_t=(e^{-ikx+\omega_jt}(i\psi^{(j)}_x-k\psi^{(j)}))_x,~~~~ x_{j-1}<x<x_{j},
\end{equation}

\noindent where $\omega_j(k)=i(\alpha_j+k^2)$ for $1\leq j\leq n+1$ and $x_0=-\infty$ and $x_{n+1}=\infty$.   Applying Green's Theorem and integrating over the (possibly unbounded) strips $(x_{j-1},x_{j})\times(0,t)$ for $1\leq j\leq n+1$, we have the $n+1$ global relations

\begin{align*}
\int_{x_{j-1}}^{x_{j}} e^{-ikx+\omega_jt}\psi^{(j)}(x,t)\ud x=&\int_{x_{j-1}}^{x_{j}} e^{-ikx}\psi^{(j)}_0(x)\ud x+\int_0^t e^{-ikx_j+ \omega_j s}(i\psi^{(j)}_x(x_j,s)-k\psi^{(j)}(x_j,s))\ud s\\
&-\int_0^t e^{-ikx_{j-1} +\omega_j s}(i\psi^{(j)}_x(x_{j-1},s)-k\psi^{(j)}(x_{j-1},s))\ud s.
\end{align*}

As before, we define the following for $j=1, \ldots, n+1$:

\begin{align*}
&\hat{\psi}^{(j)}(k,t)=\int_{x_{j-1}}^{x_j} e^{-ikx}\psi^{(j)}(x,t)\ud x, &&x_{j-1}<x<x_j, &t>0,\\
&\hat{\psi}^{(j)}_0(k)=\int_{x_{j-1}}^{x_j} e^{-ikx}\psi^{(j)}_0(x)\ud x, &&x_{j-1}<x<x_j,\\
&g_0^{(j)}(\omega,t)=\int_0^t e^{\omega s} \psi^{(j)}(x_j,s)\ud s=\int_0^t e^{\omega s} \psi^{(j+1)}(x_j,s)\ud s, &&&t>0\\
&g_1^{(j)}(\omega,t)=\int_0^t e^{\omega s} \psi^{(j)}_x(x_j,s)\ud s=\int_0^t e^{\omega s} \psi^{(j+1)}_x(x_j,s)\ud s, &&&t>0.
\end{align*}

\noindent For convenience we assume the interfaces are shifted such that $x_1=0$ and $x_j>0$ for all $j\geq2$. All but four of these integrals are proper integrals, and they are defined for $k\in \mathbb{C}$. The only ones that are not valid in all of $\mathbb{C}$ are $\hat{\psi}^{(1)}(k,t)$, $\hat{\psi}_0^{(1)}(k)$ (valid for $\Im(k)\geq 0$) and $\hat{\psi}^{(n+1)}(k,t)$, $\hat{\psi}_0^{(n+1)}(k)$ (valid for $\Im(k)\leq 0$).

With these definitions the global relations become

\begin{subequations}\label{GR_n}
\begin{align}
e^{\omega_1 t}\hat{\psi}^{(1)}(k,t)=&\hat{\psi}^{(1)}_0(k)+ig^{(1)}_1(\omega_1,t)-kg^{(1)}_0(\omega_1,t),
& \Im(k) \geq 0,\\
e^{\omega_j t}\hat{\psi}^{(j)}(k,t)=&\hat{\psi}^{(j)}_0(k)+e^{-ikx_{j}}(ig^{(j)}_1(\omega_j,t)-kg^{(j)}_0(\omega_j,t))\nonumber\\
&-e^{-ikx_{j-1} }(ig^{(j-1)}_1(\omega_j,t)-kg^{(j-1)}_0(\omega_j,t)), &k\in\CC,\\
e^{\omega_{n+1} t}\hat{\psi}^{(n+1)}(k,t)=&\hat{\psi}^{(n+1)}_0(k)-e^{-ikx_{n}}(ig^{(n)}_1(\omega_{n+1},t)-kg^{(n)}_0(\omega_{n+1},t)),
& \Im(k) \leq 0,
\end{align}
\end{subequations}
where $2\leq j\leq n$.  As in the previous section we transform the global relations so that $g^{(j)}_0(\cdot,t)$ and $g^{(j)}_1(\cdot,t)$ depend on a common argument. Let

$$
\nu^{(j)}(k)= i k\sqrt{1+\frac{\alpha_j}{k^2}},~~~~j=1, \ldots, n+1.
$$

\noindent Using the transformations $k=\pm\nu^{(j)}(\kappa)$, we have the transformed global relations

\begin{subequations}\label{GRt_n}
\begin{align}
e^{-i\kappa^2 t}\hat{\psi}^{(1)}(\nu^{(1)},t)=&\hat{\psi}^{(1)}_0(\nu^{(1)})+ig^{(1)}_1-\nu^{(1)}g^{(1)}_0, \label{GRt_n_1p}\\
e^{-i\kappa^2 t}\hat{\psi}^{(1)}(-\nu^{(1)},t)=&\hat{\psi}^{(1)}_0(-\nu^{(1)})+ig^{(1)}_1+\nu^{(1)}g^{(1)}_0, \label{GRt_n_1n} \\
\begin{split}
e^{-i\kappa^2 t}\hat{\psi}^{(j)}(\nu^{(j)},t)=&\hat{\psi}^{(j)}_0(\nu^{(j)})+e^{-i\nu^{(j)}x_j}(ig^{(j)}_1-\nu^{(j)}g^{(j)}_0)\\
&-e^{-i\nu^{(j)}x_{j-1}}(ig^{(j-1)}_1-\nu^{(j)}g^{(j-1)}_0), \label{GRt_n_jp}
\end{split}\\
\begin{split}
e^{-i\kappa^2 t}\hat{\psi}^{(j)}(-\nu^{(j)},t)=&\hat{\psi}^{(j)}_0(-\nu^{(j)})+e^{i\nu^{(j)}x_j}(ig^{(j)}_1+\nu^{(j)}g^{(j)}_0)\\
&-e^{i\nu^{(j)}x_{j-1}}(ig^{(j-1)}_1+\nu^{(j)}g^{(j-1)}_0), \label{GRt_n_jn}
\end{split}\\
e^{-i\kappa^2 t}\hat{\psi}^{(n+1)}(\nu^{(n+1)},t)=&\hat{\psi}^{(n+1)}_0(\nu^{(n+1)})-e^{-i\nu^{(n+1)}x_n}(ig^{(n)}_1-\nu^{(n+1)}g^{(n)}_0),\label{GRt_n_np}\\
e^{-i\kappa^2 t}\hat{\psi}^{(n+1)}(-\nu^{(n+1)},t)=&\hat{\psi}^{(n+1)}_0(-\nu^{(n+1)})-e^{i\nu^{(n+1)}x_n}(ig^{(n)}_1+\nu^{(n+1)}g^{(n)}_0), \label{GRt_n_nn}
\end{align}
\end{subequations}

\noindent where $2\leq j\leq n$ and $g_0^{(j)}=g_0^{(j)}(-i\kappa^2,t)$, $g_1^{(j)}=g_1^{(j)}(-i\kappa^2,t)$, $\nu^{(j)}=\nu^{(j)}(\kappa)$, for $1\leq j\leq n$.  In order for~\eqref{GRt_n} to be well defined $\Re(\kappa )\geq 0$ for~\eqref{GRt_n_1p} and~\eqref{GRt_n_nn}.  Similarly, $\Re(\kappa )\leq 0$ in~\eqref{GRt_n_1n} and~\eqref{GRt_n_np}.  Equations~\eqref{GRt_n_jp} and~\eqref{GRt_n_jn} are valid for all $\kappa \in\CC$.

Inverting the Fourier transform in~\eqref{GR_n} we have the solution formulae

\begin{subequations}
\begin{align}
\psi^{(1)}(x,t)=&\frac{1}{2\pi}\int_{-\infty}^\infty e^{ikx-\omega_1t}\hat{\psi}^{(1)}_0(k)\ud k+\frac{1}{2\pi}\int_{-\infty}^\infty e^{ikx-\omega_1t} \left(ig^{(1)}_1(\omega_1,t)-kg^{(1)}_0(\omega_1,t)\right)\ud k,\\
\begin{split}\label{middle}
\psi^{(j)}(x,t)=&\frac{1}{2\pi}\int_{-\infty}^\infty \!\!\!e^{ikx-\omega_j t}\hat{\psi}^{(j)}_0(k)\ud k +\frac{1}{2\pi}\int_{-\infty}^\infty \!\!\!\!e^{ik(x-x_j)-\omega_jt} \left(ig_1^{(j)}(\omega_j,t)-kg_0^{(j)}(\omega_j,t)\right)\!\!\ud k\\
&-\frac{1}{2\pi}\int_{-\infty}^\infty e^{ik(x-x_{j-1})-\omega_jt} \left(ig_1^{(j-1)}(\omega_j,t)-kg_0^{(j-1)}(\omega_j,t)\right)\ud k,
\end{split}\\
\begin{split}
\psi^{(n+1)}(x,t)=&\frac{1}{2\pi}\int_{-\infty}^\infty e^{ikx-\omega_{n+1}t}\hat{\psi}^{(n+1)}_0(k)\ud k\\
&-\frac{1}{2\pi}\int_{-\infty}^\infty e^{ik(x-x_n)-\omega_{n+1}t}\left(i g^{(n)}_1(\omega_{n+1},t)-kg^{(n)}_0(\omega_{n+1},t)\right)\ud k,
\end{split}
\end{align}
\end{subequations}

\noindent for $2\leq j\leq n$, and $x_{j-1}<x<x_j$.  We deform these integrals into the complex plane.  Using Cauchy's Theorem and Jordan's Lemma we have

\begin{subequations}
\begin{align}
\psi^{(1)}(x,t)=&\frac{1}{2\pi}\int_{-\infty}^\infty e^{ikx-\omega_1t}\hat{\psi}^{(1)}_0(k)\ud k-\frac{1}{2\pi}\int_{\partial D_R^{(3)}} e^{ikx-\omega_1t} \left(ig^{(1)}_1(\omega_1,t)-kg^{(1)}_0(\omega_1,t)\right)\ud k,\\
\begin{split}
\psi^{(j)}(x,t)=&\frac{1}{2\pi}\int_{-\infty}^\infty \!\!\!e^{ikx-\omega_j t}\hat{\psi}^{(j)}_0(k)\ud k
-\frac{1}{2\pi}\int_{\partial D_R^{(3)}} \!\!\!e^{ik(x-x_j)-\omega_jt} \left(ig_1^{(j)}(\omega_j,t)-kg_0^{(j)}(\omega_j,t)\right)\!\!\ud k\\
&-\frac{1}{2\pi}\int_{\partial D_R^{(1)}} e^{ik(x-x_{j-1})-\omega_jt} \left(ig_1^{(j-1)}(\omega_j,t)-kg_0^{(j-1)}(\omega_j,t)\right)\ud k,
\end{split}\\
\begin{split}
\psi^{(n+1)}(x,t)=&\frac{1}{2\pi}\int_{-\infty}^\infty e^{ikx-\omega_{n+1}t}\hat{\psi}^{(n+1)}_0(k)\ud k\\
&-\frac{1}{2\pi}\int_{\partial D_R^{(1)}} e^{ik(x-x_n)-\omega_{n+1}t}\left(i g^{(n)}_1(\omega_{n+1},t)-kg^{(n)}_0(\omega_{n+1},t)\right)\ud k,
\end{split}
\end{align}
\end{subequations}
where $D_R^{(j)}$ is as in~\eqref{DjR} and Figure~\ref{fig:LSp_DRpm}.  Again, we wish to transform the integrals involving $g_0^{(j)}(\cdot,t)$ and $g_1^{(j)}(\cdot,t)$ in each of the solution formulae above so these terms depend on $-i\kappa^2$.  As before, we deform to $D_R^{(4)}$ (with $\Lambda=\max_j{|\alpha_j|}$, $R>\sqrt{2\Lambda}$).  Choosing $k=\nu^{(j)}(\kappa)$ on $\partial D_R^{(1)}$ and $k=-\nu^{(j)}(\kappa)$ on $\partial D_R^{(3)}$ we have

\begin{subequations}\label{dsolns_n}
\begin{equation}\label{dsolns1_n}
\begin{split}
\psi^{(1)}(x,t)=&\frac{1}{2\pi}\int_{-\infty}^\infty e^{ikx-\omega_1t}\hat{\psi}^{(1)}_0(k)\ud k\\
&-\frac{1}{2\pi}\int_{\partial D_R^{(4)}} e^{-i\nu^{(1)}(\kappa)x+i\kappa^2t} \left(\frac{i\kappa}{\nu^{(1)}(\kappa)}g^{(1)}_1+\kappa g^{(1)}_0\right)\ud \kappa,
\end{split}
\end{equation}
\begin{equation}\label{dsolnsm_n}
\begin{split}
\psi^{(j)}(x,t)=&\frac{1}{2\pi}\int_{-\infty}^\infty e^{ikx-\omega_j t}\hat{\psi}^{(j)}_0(k)\ud k\\
&-\frac{1}{2\pi}\int_{\partial D_R^{(4)}} e^{-i\nu^{(j)}(\kappa)(x-x_j)+i\kappa^2t} \left(\frac{i\kappa}{\nu^{(j)}(\kappa)}g_1^{(j)}+\kappa g_0^{(j)}\right)\ud \kappa\\
&+\frac{1}{2\pi}\int_{\partial D_R^{(4)}} e^{i\nu^{(j)}(\kappa)(x-x_{j-1})+i\kappa^2 t} \left(\frac{i\kappa}{\nu^{(j)}(\kappa)}g_1^{(j-1)}-\kappa g_0^{(j-1)}\right)\ud \kappa,
\end{split}
\end{equation}
\begin{equation}\label{dsolnsn_n}
\begin{split}
\psi^{(n+1)}(x,t)=&\frac{1}{2\pi}\int_{-\infty}^\infty e^{ikx-\omega_{n+1}t}\hat{\psi}^{(n+1)}_0(k)\ud k\\
&+\frac{1}{2\pi}\int_{\partial D_R^{(4)}} e^{i\nu^{(n+1)}(\kappa)(x-x_n)+i\kappa^2t}\left(\frac{i\kappa}{\nu^{(n+1)}(\kappa)}g^{(n)}_1-\kappa g^{(n)}_0\right)\ud \kappa,
\end{split}
\end{equation}
\end{subequations}
where $g_0^{(j)}\equiv g_0^{(j)}(-i\kappa^2,t)$ and $g_1^{(j)}\equiv g_1^{(j)}(-i\kappa^2,t)$.

Using the $2n$ transformed global relations valid in $D_R^{(4)}$~\eqref{GRt_n_1p}, \eqref{GRt_n_jp}, \eqref{GRt_n_jn}, and \eqref{GRt_n_nn} one solves for $g_0^{(j)}$ and $g_1^{(j)}$. This amounts to solving the $2n\times2n$ matrix problem
\begin{equation}\label{linprob}
\mathcal{A}(\kappa)X(-i\kappa^2,t)=Y(\kappa)+e^{-i\kappa^2t}\mathcal{Y}(\kappa,t),
\end{equation}
where
\begin{subequations}\label{XYA}
\begin{equation}\label{Xvector}
X(-i\kappa^2,t)=\left(g_0^{(1)},g_0^{(2)},\ldots,g_0^{(n)},ig_1^{(1)},ig_1^{(2)},\ldots,ig_1^{(n)}\right)^\top,
\end{equation}
\begin{equation}\label{Yvector}
Y(\kappa)=-\left(\hat{\psi}_0^{(1)}(\nu^{(1)}),\ldots,\hat{\psi}_0^{(n)}(\nu^{(n)}),\hat{\psi}_0^{(2)}(-\nu^{(2)}),\ldots,\hat{\psi}_0^{(n+1)}(-\nu^{(n+1)})\right)^\top,
\end{equation}
\begin{equation}
\mathcal{Y}(\kappa,t)=\left(\hat{\psi}^{(1)}(\nu^{(1)},t),\ldots,\hat{\psi}^{(n)}(\nu^{(n)},t),\hat{\psi}^{(2)}(-\nu^{(2)},t),\ldots,\hat{\psi}^{(n+1)}(-\nu^{(n+1)},t)\right)^\top,
\end{equation}
and

\begin{equation}\label{Amatrix}
\begin{split}
&\mathcal{A}(\kappa)=\\
&\left(
\begin{array}{ccc:ccc}
-\nu^{(1)}e^{-i\nu^{(1)}x_1}&&&e^{-i\nu^{(1)}x_1}\\
\nu^{(2)}e^{-i\nu^{(2)}x_1}&-\nu^{(2)}e^{-i\nu^{(2)}x_2}&&-e^{-i\nu^{(2)}x_1}&e^{-i\nu^{(2)}x_2}\\
\hspace{.7in}\ddots&\hspace{.7in}\ddots&&\hspace{.7in}\ddots&\hspace{.7in}\ddots\\
&\nu^{(n)}e^{-i\nu^{(n)}x_{n-1}}&-\nu^{(n)}e^{-i\nu^{(n)}x_n}&&-e^{-i\nu^{(n)}x_{n-1}}&e^{-i\nu^{(n)}x_n}\\
\hdashline
-\nu^{(2)}e^{i\nu^{(2)}x_1}&\nu^{(2)}e^{i\nu^{(2)}x_2}&&-e^{i\nu^{(2)}x_1}&e^{i\nu^{(2)}x_2}\\
\hspace{.7in}\ddots&\hspace{.7in}\ddots&&\hspace{.7in}\ddots&\hspace{.7in}\ddots\\
&-\nu^{(n)}e^{i\nu^{(n)}x_{n-1}}&\nu^{(n)}e^{i\nu^{(n)}x_n}&&-e^{i\nu^{(n)}x_{n-1}}&e^{i\nu^{(n)}x_n}\\
&&-\nu^{(n+1)}e^{i\nu^{(n+1)}x_n}&&&-e^{i\nu^{(n+1)}x_n}
\end{array}
\right),
\end{split}
\end{equation}
\end{subequations}
where all $\nu^{(j)}$ are evaluated at $\kappa$.  The matrix $\mathcal{A}(\kappa)$ is made up of four $n\times n$ blocks as indicated by the dashed lines.  The two blocks in the upper half of $\mathcal{A}(\kappa)$ are zero except for entries on the main and $-1$ diagonals.  The lower two blocks of $\mathcal{A}(\kappa)$ are zero except on the main and $+1$ diagonals.

Every term in the linear equation $\mathcal{A}(\kappa)X(-i\kappa^2,t)=Y(\kappa)$ is known.  By substituting the solutions of this equation into~\eqref{dsolns_n}, we have solved the LS equation with a piecewise constant potential in terms of only known functions.  It remains to show that the contribution to the solution from the linear equation $\mathcal{A}(\kappa)X(-i\kappa^2,t)=e^{-i\kappa^2t}\mathcal{Y}(\kappa,t)$ is 0 when substituted into~\eqref{dsolns_n}.

To this end consider $\mathcal{A}(\kappa)X(-i\kappa^2,t)=e^{-i\kappa^2t}\mathcal{Y}(\kappa,t).$  We solve this system using Cramer's Rule.  We factor $\mathcal{A}(\kappa)=\mathcal{A}^L(\kappa)\mathcal{A}^M(\kappa)$ where
\begin{equation}\label{ALmatrix}
\mathcal{A}^L(\kappa)=\left(\begin{array}{ccc:ccc}
e^{-i\nu^{(1)}x_1}&&&\\
&\ddots\\
&&e^{-i\nu^{(n)}x_n}&\\
\hdashline
&&&e^{i\nu^{(1)}x_1}\\
&&&&\ddots\\
&&&&&e^{i\nu^{(n)}x_n}
\end{array}\right),
\end{equation}
Let $\mathcal{A}_j(\kappa,t)$ be the matrix $\mathcal{A}(\kappa)$ with the $j^{\textrm{th}}$ column replaced by $e^{-i\kappa^2t}\mathcal{Y}(\kappa,t)$.  Similar to $\mathcal{A}(\kappa)$, this matrix can be factored as $\mathcal{A}_j(\kappa,t)=e^{-i\kappa^2t}\mathcal{A}^L(\kappa)\mathcal{A}_j^M(\kappa,t)$.  Hence, $\det(\mathcal{A}_j(\kappa,t))=e^{-i\kappa^2t}\det(\mathcal{A}_j^M(\kappa,t))$.

The terms we are trying to eliminate contribute to the solution~\eqref{dsolns_n} in the form
\begin{subequations}\label{badparts}
\begin{align}
&\int_{\partial D_R^{(4)}} e^{-i\nu^{(j)}(x-x_j)+i\kappa^2t} \left(\frac{i\kappa}{\nu^{(j)}}g_1^{(j)}+\kappa g_0^{(j)}\right)\ud \kappa,
\end{align}
and
\begin{align}
&\int_{\partial D_R^{(4)}} e^{i\nu^{(j)}(x-x_{j-1})+i\kappa^2 t} \left(\frac{i\kappa}{\nu^{(j)}}g_1^{(j-1)}-\kappa g_0^{(j-1)}\right)\ud \kappa ,
\end{align}
\end{subequations}
for $1\leq j\leq n$, $2\leq j\leq n+1$, respectively with $x_{j-1}<x<x_j$.  Using Cramer's Rule these become
\begin{subequations}\label{badparts_det}
\begin{align}
&\int_{\partial D_R^{(4)}} e^{-i\nu^{(j)}(x-x_j)} \left(\frac{\kappa}{\nu^{(j)}}\frac{\det(\mathcal{A}_{j+n}^M(\kappa,t))}{\det(\mathcal{A}^M(\kappa))}+\kappa  \frac{\det(\mathcal{A}_{j}^M(\kappa))}{\det(\mathcal{A}^M(\kappa))}\right)\ud \kappa ,
\end{align}
and
\begin{align}
&\int_{\partial D_R^{(4)}} e^{i\nu^{(j)}(x-x_{j-1})} \left(\frac{\kappa}{\nu^{(j)}}\frac{\det(\mathcal{A}_{j-1+n}^M(\kappa,t))}{\det(\mathcal{A}^M(\kappa))}-\kappa \frac{\det(\mathcal{A}_{j-1}^M(\kappa,t))}{\det(\mathcal{A}^M(\kappa))}\right)\ud \kappa ,
\end{align}
\end{subequations}
respectively.  Here we have used the factorizations of $\det(\mathcal{A}(\kappa))$ and $\det(\mathcal{A}_j(\kappa,t))$.   As is usual in the UTM we use the large $\kappa $ asymptotics to show the terms~\eqref{badparts_det} are 0.

Observe the elements of $\mathcal{A}^M(\kappa)$ are either 0, $\mathcal{O}(\kappa)$ or decaying exponentially fast for $\mathcal{\kappa}\in D_R^{(4)}$ and all but $n$ columns are $\mathcal{O}(1)$.  Hence,
$$\det(\mathcal{A}^M(\kappa))=\det(\mathcal{A}(\kappa))=\Delta(\kappa)=\mathcal{O}(\kappa^n),$$ for large $\kappa$ in $D_R^{(4)}$.  Expanding the determinant of $\mathcal{A}^M_j(\kappa,t)$ along the $j^{\textrm{th}}$ column we see that

\begin{align*}
e^{-i\nu^{(j)}(x-x_j)}\kappa \frac{\det(\mathcal{A}_{j}^M(\kappa,t))}{\det(\mathcal{A}^M(\kappa))}=& e^{-i\nu^{(j)}(x-x_j)}\kappa \frac{\det(\mathcal{A}_{j}^M(\kappa,t))}{\Delta(\kappa)}\\
=& e^{-i\nu^{(j)}(x-x_j)}\sum_{\ell=1}^n c_{\ell}(\kappa)e^{ix_\ell\nu^{(\ell)}}\hat{\psi}^{(\ell)}\left(\nu^{(\ell)},t\right)\\
&\hspace{1in}+c_{\ell+n}(\kappa)e^{-ix_\ell\nu^{(\ell)}}\hat{\psi}^{(\ell+1)}\left(-\nu^{(\ell+1)},t\right),
\end{align*}
where $c_{\ell}(\kappa)=\mathcal{O}(\kappa^{0})$ and $x_{j-1}<x<x_j$. 
In the large $|\kappa|$ limit, $\nu^{(j)}(\kappa)\sim -i\kappa$, so $\mathcal{A}^M(\kappa)$ reduces to the value of $\mathcal{A}^M(\kappa)$ with $\alpha=0$.
Note $e^{ix_\ell\nu^{(\ell)}}\hat{\psi}^{(\ell)}\left(\nu^{(\ell)},t\right)$ and $e^{-ix_\ell\nu^{(\ell)}}\hat{\psi}^{(\ell+1)}\left(-\nu^{(\ell+1)},t\right)$ decay uniformly for $\kappa$ in the closure of $D_R^{(4)}$. The integrands in~\eqref{badparts_det} are analytic for $\kappa\in D_R^{(4)}$. The zeros of $\det(\mathcal{A}^M(\kappa))$ are confined to strips of asymptotically constant width that are parallel to either the real or imaginary axis (depending on the sign of $\alpha_j$)~\cite{Langer1929, Langer}.  
In the examples we consider, by choosing a sufficiently large $R$ one is able to choose a region $D_R^{(4)}$ where~\eqref{badparts_det} are analytic.  The analysis of the zeros of $\Delta(\kappa)$ in its full generality is difficult and is not attempted here.
Similar to the argument on page~\pageref{JordanCauchyArgument}, since the integral along $\mathcal{L}^{(4)}_C$ vanishes for large $C$, the integrals~\eqref{badparts_det} must vanish since the contour $\mathcal{L}_{D^{(4)}}$  becomes $\partial D^{(4)}$ as $C\to\infty$.  The uniform decay of the ratios of the determinants for large $\kappa$ is exactly the condition required for the integral to vanish using Jordan's Lemma.  Hence, the solution to~\eqref{ls_p_n} is~\eqref{dsolns_n} where $g_0^{(j)}(-i\kappa^2,t)$ and $g_1^{(j)}(-i\kappa^2,t)$ for $1\leq j\leq n+1$ are found by solving
\begin{equation}\label{linear_system}
\mathcal{A}(\kappa)X(-i\kappa^2,t)=Y(\kappa),
\end{equation} where $\mathcal{A}(\kappa), X(-i\kappa^2,t)$, and $Y(\kappa)$ are given in Equations~\eqref{Amatrix},~\eqref{Xvector}, and~\eqref{Yvector} respectively.  As in the previous section, deforming to the real line is possible using~\eqref{1_imag}-\eqref{3_real}.  However, one must be careful to avoid any poles present in~\eqref{dsolns_n}.


\section{Potential well and barrier}\label{sec:well}
As an example of the general method given in Section~\ref{sec:n}, in this section we solve the classical problem of the finite potential well or barrier:
\begin{equation}\label{ls_p_well}
i\psi_t=-\psi_{xx}+\alpha(x)\psi,
\end{equation}
for $-\infty<x<\infty$ and $$\alpha(x)=\left\{\begin{array}{lcrr} 0,&&&x<x_1,\\\alpha,&&x_1<&x<x_2,\\ 0,&&&x>x_2, \end{array}\right.$$ with the initial condition $\psi(x,0)=\psi_0(x)$ and $\lim_{|x|\to\infty}\psi(x,t)=0$ with $\psi(x,t)$ and its spatial derivative in $L^1$.

The problem of a finite potential well or barrier is a standard textbook problem in quantum mechanics~\cite{FeynmanLecture2, LandauLifshitz, Merzbacher}.  In such texts this problem is usually solved using separation of variables, \emph{i.e.}, assuming $\psi(x,t)=X(x)T(t)$.  The $x$ problem, $X''+(\xi^2-\alpha(x)) X=0$ is solved in the three different regions.  Separation of variables is only useful if the initial wave function $\psi(x,0)$ can be expanded in terms of solutions of the time-independent Schr\"odinger equation~\cite{Merzbacher}.  Solving the time-independent Schr\"odinger equation is equivalent to studying the forward scattering problem with the specified potential.  The ``scattering matrix"  (see \cite[Equation 1.3.3]{AS} or \cite[p.\ 104]{Faddeyev}) is $$\left(\begin{array}{lr} a(\xi)&b(\xi)\\\overline{b}(\xi) &-\overline{a}(\xi)\end{array}\right).$$  The zeros of $a(\xi)$ are the discrete eigenvalues for the problem. With some straightforward work we find \begin{equation}\label{aeqn}a(\xi )=e^{i\xi  x_2}\left(\cosh(x_2\sqrt{\alpha-\xi^2})-\frac{i(2\xi^2-\alpha)}{2\xi \sqrt{\alpha-\xi^2}}\sinh(x_2\sqrt{\alpha-\xi^2})\right).\end{equation} This problem is examined in many excellent texts including~\cite{AS, AgranovichMarchenko, DeiftTrubowitz, Faddeyev}.

The potential well or barrier problem is the standard example to introduce students to the concept of quantum tunneling which is a phenomenon where a particle ``tunnels" over a barrier that it cannot overcome in the classical mechanics setting~\cite{Razavy}. The closed form solutions we present at the end of this section all depend on the initial conditions from each of the three regions and quantum tunneling is clearly present.

\begin{figure}
\begin{center}
\def\svgwidth{\textwidth}
   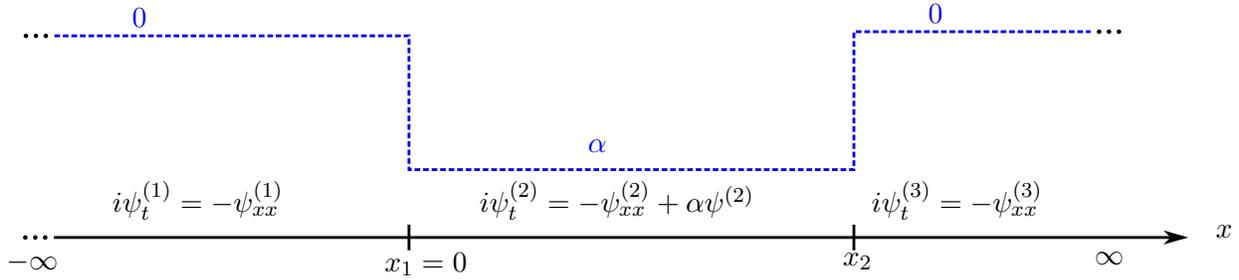 
   \caption{A cartoon of the potential $\alpha(x)$ for a potential well or barrier.   \label{fig:wellbarrier}}
   \end{center}
\end{figure}
Finding the closed form solutions is as easy as letting $n=2$, $\alpha_1=\alpha_3=0$ and $\alpha_2=\alpha$ in~\eqref{dsolns_n} as in Figure~\ref{fig:wellbarrier}.  Again we denote $g_0^{(j)}=g_0^{(j)}(-i\kappa^2,t)$, $g_1^{(j)}=g_1^{(j)}(-i\kappa^2,t)$, for $j=1,2$ and $\nu^{(j)}=\nu^{(j)}(\kappa)$.  Solving~\eqref{linear_system} in the case of $n=2$ we have solutions for $g_0^{(1)}, g_1^{(1)}, g_0^{(2)}, g_1^{(2)}$ valid in $D_R^{(4)}$. Let
\begin{align*}
\Delta(\kappa)=&2\pi \left(i\kappa(e^{-ix_2\nu^{(2)}}+1)+\nu^{(2)}(e^{-ix_2\nu^{(2)}}-1)\right)\left(i\kappa(1-e^{ix_2\nu^{(2)}})+\nu^{(2)}(1+e^{ix_2\nu^{(2)}})\right)\\
=&4i\pi \left( (\alpha+2\kappa^2)\sin(x_2\nu^{(2)})+2\kappa\nu^{(2)}\cos(x_2\nu^{(2)})\right).
\end{align*}
The solutions~\eqref{dsolns_n} with the appropriate values of $g_0^{(j)}$ and $g_1^{(j)}$ are

\begin{align}
\begin{split}\label{fullsolns1_3}
\psi^{(1)}(x,t)=&\frac{1}{2\pi}\int_{-\infty}^\infty e^{ikx-\omega_1t}\hat{\psi}^{(1)}_0(k)\ud k+\int_{\partial D_R^{(4)}}  \frac{i\alpha(e^{-ix_2\nu^{(2)}}-e^{ix_2\nu^{(2)}})}{\Delta(\kappa)} e^{\kappa x+i\kappa^2t}  \hat{\psi}^{(1)}_0(i\kappa ) \ud \kappa \\
&+\int_{\partial D_R^{(4)}}  \frac{2i\kappa(\kappa-i\nu^{(2)}) }{\Delta(\kappa )} e^{\kappa x-ix_2\nu^{(2)}+i\kappa^2t} \hat{\psi}^{(2)}_0(-\nu^{(2)}) \ud \kappa\\
&+\int_{\partial D_R^{(4)}}  \frac{2i\kappa (\kappa+i\nu^{(2)})}{\Delta(\kappa )}e^{\kappa (x2ix_2\nu^{(2)})+i\kappa^2t} \hat{\psi}^{(2)}_0(\nu^{(2)}) \ud \kappa \\
&-\int_{\partial D_R^{(4)}}  \frac{4\kappa\nu^{(2)} }{\Delta(\kappa )}e^{\kappa (x+x_2)+i\kappa^2t} \hat{\psi}^{(3)}_0(-i\kappa ) \ud \kappa ,
\end{split}
\end{align}
for $x<0$,

\begin{align}
\begin{split}\label{fullsolns2_3}
\psi^{(2)}(x,t)=&\frac{1}{2\pi}\int_{-\infty}^\infty e^{ikx-\omega_2 t}\hat{\psi}^{(2)}_0(k)\ud k+\int_{\partial D_R^{(4)}}\frac{2i\kappa(\kappa+i\nu^{(2)})}{\Delta(\kappa )} e^{i\nu^{(2)}(x_2-x)+i\kappa^2t} \hat{\psi}^{(1)}_0(i\kappa ) \ud \kappa \\
&-\int_{\partial D_R^{(4)}}\frac{\kappa(\kappa+i\nu^{(2)})^2}{\nu^{(2)}\Delta(\kappa )} e^{i\nu^{(2)}(x_2-x)+i\kappa^2t}  \hat{\psi}^{(2)}_0(-\nu^{(2)}) \ud \kappa \\
&+\int_{\partial D_R^{(4)}}\frac{\alpha\kappa}{\nu^{(2)}\Delta(\kappa )}e^{i\nu^{(2)}(x_2-x)+i\kappa^2t}  \hat{\psi}^{(2)}_0(\nu^{(2)}) \ud \kappa \\
&-\int_{\partial D_R^{(4)}}\frac{2i\kappa(\kappa-i\nu^{(2)})}{\Delta(\kappa )}e^{-i\nu^{(2)}x+\kappa x_2+i\kappa^2t} \hat{\psi}^{(3)}_0(-i\kappa ) \ud \kappa \\
&-\int_{\partial D_R^{(4)}}\frac{2i\kappa(\kappa-i\nu^{(2)})}{\Delta(\kappa )} e^{i\nu^{(2)}(x-x_2)+i\kappa^2t} \hat{\psi}^{(1)}_0(i\kappa ) \ud \kappa \\
&+\int_{\partial D_R^{(4)}}\frac{\alpha \kappa}{\nu^{(2)}\Delta(\kappa )}e^{i\nu^{(2)}(x-x_2)+i\kappa^2t}  \hat{\psi}^{(2)}_0(-\nu^{(2)}) \ud \kappa \\
&+\int_{\partial D_R^{(4)}}\frac{\kappa(\kappa+i\nu^{(2)})^2}{\nu^{(2)}\Delta(\kappa )} e^{i\nu^{(2)}(x_2+x)+i\kappa^2t}\hat{\psi}^{(2)}_0(\nu^{(2)}(\kappa )) \ud \kappa \\
&+\int_{\partial D_R^{(4)}}\frac{2i\kappa (\kappa+i\nu^{(2)})}{\Delta(\kappa )}e^{i\nu^{(2)}x+\kappa x_2+i\kappa^2t} \hat{\psi}^{(3)}_0(-i\kappa ) \ud \kappa ,\end{split}
\end{align}
for $0<x<x_2$, and
\begin{align}
\begin{split}\label{fullsolns3_3}
\psi^{(3)}(x,t)=&\frac{1}{2\pi}\int_{-\infty}^\infty e^{ikx-\omega_{3}t}\hat{\psi}^{(3)}_0(k)\ud k-\int_{\partial D_R^{(4)}}\frac{4\kappa \nu^{(2)} }{\Delta(\kappa )}e^{\kappa (x_2-x)+i\kappa^2t} \hat{\psi}^{(1)}_0(i\kappa ) \ud \kappa \\
&+\int_{\partial D_R^{(4)}}\frac{2\kappa(\kappa+i\nu^{(2)}) }{\Delta(\kappa )}e^{\kappa (x_2-x)+i\kappa^2t}\hat{\psi}^{(2)}_0(-\nu^{(2)}) \ud \kappa \\
&-\int_{\partial D_R^{(4)}}\frac{2i\kappa(\kappa-i\nu^{(2)}) }{\Delta(\kappa )} e^{\kappa (x_2-x)+i\kappa^2t} \hat{\psi}^{(2)}_0(\nu^{(2)}(\kappa )) \ud \kappa \\
&+\int_{\partial D_R^{(4)}}\frac{i\alpha\left(1-e^{2i\nu^{(2)}x_2}\right)}{\Delta(\kappa )}e^{\kappa (x_2-x)+i\kappa^2t} \hat{\psi}^{(3)}_0(-i\kappa ) \ud \kappa,
\end{split}
\end{align}
when $x>x_2$.

Using the change of variables $\kappa=ik$ in~\eqref{fullsolns1_3}, $\kappa=-ik$ in~\eqref{fullsolns3_3}, $\kappa=ik\sqrt{1+\frac{\alpha}{k^2}}$ in the second, third, fourth, and fifth integrals of~\eqref{fullsolns2_3}, and $\kappa=-ik\sqrt{1+\frac{\alpha}{k^2}}$ in the last four integrals of~\eqref{fullsolns2_3} we find
\begin{align}
\begin{split}\label{t_fullsolns1_3}
\psi^{(1)}(x,t)=&\frac{1}{2\pi}\int_{-\infty}^\infty e^{ikx-\omega_1t}\hat{\psi}^{(1)}_0(k)\ud k\\
&+\int_{\partial D_R^{(3)} }  \frac{\alpha \left(e^{-ikx_2\sqrt{1-\frac{\alpha}{k^2}}}-e^{ikx_2\sqrt{1-\frac{\alpha}{k^2}}} \right) }{\Delta(ik)}e^{ik x-\omega_1 t}\hat{\psi}^{(1)}_0(-k) \ud k\\
&-\int_{\partial D_R^{(3)} }  \frac{2k^2\left(1+\sqrt{1-\frac{\alpha}{k^2}}\right)}{\Delta(ik)}  e^{ik\left(x+x_2\sqrt{1-\frac{\alpha}{k^2}}\right)-\omega_1t}\hat{\psi}^{(2)}_0\left(k\sqrt{1-\frac{\alpha}{k^2}}\right) \ud k\\
&+\int_{\partial D_R^{(3)} }  \frac{2k^2\left(1-\sqrt{1-\frac{\alpha}{k^2}}\right)}{\Delta(ik)} e^{ik\left(x-x_2\sqrt{1-\frac{\alpha}{k^2}}\right)-\omega_1t}\hat{\psi}^{(2)}_0\left(-k\sqrt{1-\frac{\alpha}{k^2}}\right)  \ud k\\
&-\int_{\partial D_R^{(3)} }  \frac{4k^2\sqrt{1-\frac{\alpha}{k^2}}  }{\Delta(ik)}  e^{ik\left(x+x_2\right)-\omega_1t}\hat{\psi}^{(3)}_0(k)\ud k,
\end{split}
\end{align}
for $x<0$,

\begin{align}
\begin{split}\label{t_fullsolns2_3}
\psi^{(2)}(x,t)=&\frac{1}{2\pi}\int_{-\infty}^\infty e^{ikx-\omega_2 t}\hat{\psi}^{(2)}_0(k)\ud k\\
&-\int_{\partial D_R^{(3)}}\frac{2k^2(1-\sqrt{1+\frac{\alpha}{k^2}})}{\Delta\left(ik\sqrt{1+\frac{\alpha}{k^2}}\right)}e^{ik (x-x_2) -\omega_2t}\hat{\psi}^{(1)}_0\left(-k\sqrt{1-\frac{\alpha}{k^2}}\right)  \ud k\\
&+\int_{\partial D_R^{(3)}}\frac{k^2(\sqrt{1+\frac{\alpha}{k^2}}-1)^2}{\Delta\left(ik\sqrt{1+\frac{\alpha}{k^2}}\right)} e^{ik(x-x_2)-\omega_2t}  \hat{\psi}^{(2)}_0(k)\ud k\\
&+\int_{\partial D_R^{(3)}}\frac{\alpha }{\Delta\left(ik\sqrt{1+\frac{\alpha}{k^2}}\right)} e^{ik(x-x_2)-\omega_2t}\hat{\psi}^{(2)}_0(-k) \ud k\\
&-\int_{\partial D_R^{(3)}}\frac{2k^2(1+\sqrt{1+\frac{\alpha}{k^2}})}{\Delta\left(ik\sqrt{1+\frac{\alpha}{k^2}}\right)}e^{ik (x+x_2\sqrt{1+\frac{\alpha}{k^2}})-\omega_2t }\hat{\psi}^{(3)}_0\left(k\sqrt{1-\frac{\alpha}{k^2}}\right)  \ud k\\
&+\int_{\partial D_R^{(1)}}\frac{2k^2(1+\sqrt{1+\frac{\alpha}{k^2}})}{\Delta\left(-ik\sqrt{1+\frac{\alpha}{k^2}}\right)} e^{ik(x-x_2)-\omega_2t}\hat{\psi}^{(1)}_0\left(k\sqrt{1-\frac{\alpha}{k^2}} \right) \ud k\\
&+\int_{\partial D_R^{(1)}}\frac{k^2(\sqrt{1+\frac{\alpha}{k^2}}-1)^2}{\Delta\left(-ik\sqrt{1+\frac{\alpha}{k^2}}\right)}e^{ik(x+x_2)-\omega_2t}\hat{\psi}^{(2)}_0(k)  \ud k\\
&-\int_{\partial D_R^{(1)}}\frac{\alpha }{\Delta\left(-ik\sqrt{1+\frac{\alpha}{k^2}}\right)}  e^{ik(x-x_2)-\omega_2t}\hat{\psi}^{(2)}_0(-k)\ud k\\
&+\int_{\partial D_R^{(1)}}\frac{2k^2(1-\sqrt{1+\frac{\alpha}{k^2}})}{\Delta\left(-ik\sqrt{1+\frac{\alpha}{k^2}}\right)} e^{ik (x-x_2\sqrt{1+\frac{\alpha}{k^2}})-\omega_2t}\hat{\psi}^{(3)}_0\left(-k\sqrt{1-\frac{\alpha}{k^2}}\right) \ud k,
\end{split}
\end{align}
for $0<x<x_2$, and
\begin{align}
\begin{split}\label{t_fullsolns3_3}
\psi^{(3)}(x,t)=&\frac{1}{2\pi}\int_{-\infty}^\infty e^{ikx-\omega_{3}t}\hat{\psi}^{(3)}_0(k)\ud k\\
&+\int_{\partial D_R^{(1)}}\frac{4k^2\sqrt{1-\frac{\alpha}{k^2}}}{\Delta(-ik)} e^{ik(x-x_2) -\omega_3 t }\hat{\psi}^{(1)}_0(k)  \ud k\\
&+\int_{\partial D_R^{(1)}}\frac{2k^2(1+\sqrt{1-\frac{\alpha}{k^2}})}{\Delta(-ik)} e^{ik(x-x_2)-\omega_3t }\hat{\psi}^{(2)}_0\left(k\sqrt{1-\frac{\alpha}{k^2}} \right) \ud k\\
&-\int_{\partial D_R^{(1)}}\frac{2k^2(1-\sqrt{1-\frac{\alpha}{k^2}})}{\Delta(-ik)} e^{ik (x-x_2)-\omega_3t} \hat{\psi}^{(2)}_0\left(-k\sqrt{1-\frac{\alpha}{k^2}}\right)  \ud k\\
&+\int_{\partial D_R^{(1)}}\frac{\alpha(e^{-ikx_2\sqrt{1-\frac{\alpha}{k^2}}}-e^{ikx_2\sqrt{1-\frac{\alpha}{k^2}}}) }{\Delta(-ik)}e^{ik(x-2x_2)-\omega_3t}\hat{\psi}^{(3)}_0(-k)  \ud k,
\end{split}
\end{align}
when $x>x_2$.

\medskip
{\bf Remarks:}
\begin{itemize}

\item If one lets $\alpha=0$ in~\eqref{t_fullsolns1_3}-\eqref{t_fullsolns3_3} then the Fourier transform solution to the free Schr\"odinger equation on the whole line is recovered.

\item In order to numerically or asymptotically evaluate these expressions one could use techniques presented in~\cite{BiondiniTrogdon, Levin, TrogdonThesis, Trogdon}.  

\item As stated at the beginning of this section,~\eqref{ls_p_well} is solved in standard quantum mechanics texts using separation of variables and the study of the forward scattering problem with the specified potential.  The zeros of $a(\xi)$, the $(1,1)$ component of the scattering matrix, are the discrete eigenvalues for the problem.   The zeros of $a(\xi)$ cannot be found explicitly but it is clear that the zeros of $a(\xi)$ for $\xi$ purely imaginary correspond to the zeros of the denominators of~\eqref{t_fullsolns1_3}-\eqref{t_fullsolns3_3} with $i\xi^2=\omega_j(k)$.  The contribution of the discrete spectrum can be recovered explicitly by deforming contours in the complex plane to the real line, resulting in a sum of residue contributions corresponding to the eigenmodes of the problem.  This is done explicitly for other problems solved via the UTM in~\cite{FokasBook}.

\item As is typical in using the UTM, we find the solution under the assumption of existence.  Often, to justify existence, one checks that the solution formula obtained actually satisfies the original problem a posteriori. This is not attempted in this paper.

\end{itemize}

\section{Initial-to-Interface Map}\label{sec:i2i}
The construction of a Dirichlet-to-Neumann map, that is, determining the boundary values that are not prescribed in terms of the initial and boundary conditions, is important in the study of PDEs and particularly in inverse problems~\cite{Fokas8, SylvesterUhlmann}.  In this section we construct a similar map between the initial values of the PDE and the function (and its first spatial derivative) evaluated at the interfaces.  This map allows for an alternative to the approach of finding solutions to interface problems as presented in the first four sections of this paper.  Given the initial conditions one could find the value of the function and its derivatives at the interface(s).  This changes the problem at hand from an interface problem to a (consistently) overspecified  BVP.  At this point, the BVP can be solved on any segment using any number of methods appropriate to the given problem.

In this section we construct the initial-to-interface map for the IVP~\eqref{ls_p_n}.  We begin by evaluating the $2n\times 2n$ linear equation~\eqref{linprob} at $t=T$:
\begin{equation}\label{linprobT}
\mathcal{A}(\kappa)X(-i\kappa^2,T)=Y(\kappa)+e^{-i\kappa^2T}\mathcal{Y}(\kappa,T),
\end{equation}
where $\mathcal{A}(\kappa)$, $X(-i\kappa^2,T)$, $Y(\kappa)$, and $\mathcal{Y}(\kappa,T)$ are given in~\eqref{XYA}.  Using Cramer's Rule to solve this system we have
\begin{subequations}\label{I2I_g01solns}
\begin{align}
g_0^{(j)}(-i\kappa^2,T)=&\frac{\det(\mathcal{A}_j(\kappa,T))}{\det(\mathcal{A}(\kappa))},\\
g_1^{(j)}(-i\kappa^2,T)=&-i\frac{\det(\mathcal{A}_{j+n}(\kappa,T))}{\det(\mathcal{A}(\kappa))},
\end{align}
\end{subequations}
where $1\leq j\leq n$ and $\mathcal{A}_j(\kappa,T)$ is the matrix $\mathcal{A}(\kappa)$ with the $j^\textrm{th}$ column replaced by $Y(\kappa)+e^{-i\kappa^2T}\mathcal{Y}(\kappa,T)$.  This does not give an effective initial-to-interface map because~\eqref{I2I_g01solns} depends on the solutions $\hat{\psi}^{(j)}(\cdot,T)$.  To eliminate this dependence we multiply~\eqref{I2I_g01solns} by $\kappa e^{i\kappa^2t}$ and integrate around $D_R^{(4)}$, as is typical in the construction of Dirichlet-to-Neumann maps~\cite{FokasBook}.  Switching the order of integration we have
\begin{subequations}\label{gsolns}
\begin{align}
\int_0^T \psi^{(j)}(x_j,s) \int_{\partial D_R^{(4)}}\kappa e^{i\kappa^2(t-s)}\ud \kappa\ud s=&\int_{\partial D_R^{(4)}} e^{i\kappa^2t} \frac{\kappa\det(\mathcal{A}_j(\kappa,T))}{\det(\mathcal{A}(\kappa))}\ud \kappa, \\
\int_0^T \psi^{(j)}_x(x_j,s) \int_{\partial D_R^{(4)}} \kappa e^{i\kappa^2(t-s)}\ud \kappa\ud s=&-i\int_{\partial D_R^{(4)}} e^{i\kappa^2t}\frac{\kappa \det(\mathcal{A}_{j+n}(\kappa,T))}{\det(\mathcal{A}(\kappa))}\ud \kappa.
\end{align}
\end{subequations}
Using the change of variables $\ell=\kappa^2$ and the classical Fourier transform formula for the delta function we find
\begin{subequations}\label{I2I_INTg01solns}
\begin{align}
\psi^{(j)}(x_j,t) =&\frac{1}{\pi}\int_{\partial D_R^{(4)}} e^{i\kappa^2t} \frac{\kappa\det(\mathcal{A}_j(\kappa,T))}{\det(\mathcal{A}(\kappa))}\ud \kappa, \\
\psi^{(j)}_x(x_j,t) =&\frac{-i}{\pi}\int_{\partial D_R^{(4)}} e^{i\kappa^2t}\frac{\kappa \det(\mathcal{A}_{j+n}(\kappa,T))}{\det(\mathcal{A}(\kappa))}\ud \kappa.
\end{align}
\end{subequations}

To examine the right-hand-side of~\eqref{I2I_INTg01solns} we factor the matrix $\mathcal{A}(\kappa)$ as $\mathcal{A}^L(\kappa)\mathcal{A}^M(\kappa)$ where $\mathcal{A}^L(\kappa)$ is given by~\eqref{ALmatrix}.  Similarly, $\mathcal{A}_j(\kappa,T)=e^{-i\kappa^2T}\mathcal{A}^L(\kappa)\mathcal{A}_j^M(\kappa,T)$.  Using these factorizations,~\eqref{I2I_INTg01solns} becomes

\begin{subequations}\label{i2i_int_solns}
\begin{align}
\psi^{(j)}(x_j,t) =&\frac{1}{\pi}\int_{\partial D_R^{(4)}} e^{i\kappa^2(t-T)} \frac{\kappa\det(\mathcal{A}_j^M(\kappa,T))}{\det(\mathcal{A}^M(\kappa))}\ud \kappa, \\
\psi^{(j)}_x(x_j,t) =&\frac{-i}{\pi}\int_{\partial D_R^{(4)}} e^{i\kappa^2(t-T)}\frac{\kappa \det(\mathcal{A}_{j+n}^M(\kappa,T))}{\det(\mathcal{A}^M(\kappa))}\ud \kappa.
\end{align}
\end{subequations}
As in Section~\ref{sec:n}, the elements of $\mathcal{A}^M(\kappa)$ are either 0, $\mathcal{O}(\kappa)$ or decaying exponentially fast for $\mathcal{\kappa}\in D_R^{(4)}$ and $$\det(\mathcal{A}^M(\kappa))=\Delta(\kappa)=\mathcal{O}(\kappa^n),$$ for large $\kappa$ in $D_R^{(4)}$.  Expanding the determinant of $\mathcal{A}^M_j(\kappa,T)$ along the $j^{\textrm{th}}$ column we see that

\begin{align*}
e^{i\kappa^2(t-T)} \frac{\kappa \det(\mathcal{A}_{j}^M(\kappa,T))}{\det(\mathcal{A}^M(\kappa))}=& e^{i\kappa^2(t-T)} \frac{\kappa}{c(\kappa)}\det(\mathcal{A}_{j}^M(\kappa,T))\\
=& e^{i\kappa^2(t-T)} \sum_{\ell=1}^n c_{\ell}(\kappa)e^{ix_\ell\nu^{(\ell)}}\hat{\psi}^{(\ell)}\left(\nu^{(\ell)},T\right)\\
&\hspace{.8in}+c_{\ell+n}(\kappa)e^{-ix_\ell\nu^{(\ell)}}\hat{\psi}^{(\ell+1)}\left(-\nu^{(\ell+1)},T\right),
\end{align*}
where $c_{\ell}(\kappa)=\mathcal{O}(\kappa^{0})$ and $x_{j-1}<x<x_j$.  The terms $e^{ix_\ell\nu^{(\ell)}}\hat{\psi}^{(\ell)}\left(\nu^{(\ell)},T\right)$ and\\ $e^{-ix_\ell\nu^{(\ell)}}\hat{\psi}^{(\ell+1)}\left(-\nu^{(\ell+1)},T\right)$ decay exponentially for $k\in D_R^{(4)}$ and the integrands of~\eqref{i2i_int_solns} are analytic for $\Re(\kappa)>0$.  As in previous sections, since the integral along $\mathcal{L}^{(4)}_C$ vanishes for large $C$, the integrals~\eqref{i2i_int_solns} must vanish since the contour $\mathcal{L}_{D^{(4)}}$  becomes $\partial D^{(4)}$ as $C\to\infty$.

Since the terms involving the elements of $\mathcal{Y}(\kappa,T)$ evaluate to zero in the solution expression the initial-to-interface map for~\eqref{ls_p_n} is
\begin{subequations}\label{i2i_solns}
\begin{align}
\psi^{(j)}(x_j,t) =&\frac{1}{\pi}\int_{\partial D_R^{(4)}} e^{i\kappa^2t} \frac{\kappa\det(A_j^M(\kappa))}{\det(\mathcal{A}^M(\kappa))}\ud \kappa, \\
\psi^{(j)}_x(x_j,t) =&\frac{-i}{\pi}\int_{\partial D_R^{(4)}} e^{i\kappa^2t}\frac{\kappa \det(A_{j+n}^M(\kappa))}{\det(\mathcal{A}^M(\kappa))}\ud \kappa,
\end{align}
\end{subequations}
where $A_j(\kappa)$ is the matrix $\mathcal{A}(\kappa)$ with the $j^\textrm{th}$ column replaced by $Y(\kappa)$ and is factored as $A_j(\kappa)=\mathcal{A}^L(\kappa)A_j^M(\kappa)$.  Equation~\eqref{i2i_solns} is an effective map between the values of the function at the interface and the given initial conditions.

\medskip
\textbf{Remark.} Since the problem is linear, one could have assumed the initial condition was zero for $x$ outside the region $x_{\ell-1}<x<x_\ell$.  Then, the map would be in terms of just $\psi_0^{(\ell)}(\cdot)$.  Summing over $1\leq \ell\leq n+1$ would give the complete map for a general initial condition.

\section{Conclusion}
In this paper we find explicit, closed-form solutions to the time-dependent linear Schr\"odinger equation with a piecewise constant potential.  Further, we construct an initial-to-interface map which allows one to change the problem from an interface problem to a BVP.  This is a classical problem with important applications.

\section*{Acknowledgments} This work was generously supported by the National Science Foundation grant number NSF-DGE-0718124 (N.E.S).  Any opinions, findings, and conclusions or recommendations expressed in this material are those of the authors and do not necessarily reflect the views of the funding sources.  The authors  wish to thank the referees for their careful reading of the manuscript, correcting some errors, as well as a number of useful suggestions.

\bibliographystyle{abbrv}
\bibliography{FullBib}

\end{document}